\documentclass[12pt,a4paper]{article}

\usepackage{amsmath}
\usepackage{bbm}
\usepackage{mathrsfs}
\usepackage{latexsym}
\usepackage{amssymb}
\usepackage{graphicx}
\usepackage{verbatim}
\usepackage{psfrag}
\usepackage{sirlarticle}

\newcommand{\Scheffe}{Scheff\'{e}}

\renewcommand{\b}[1]{\left[{#1}\right]}

\renewcommand{\epsilon}{\varepsilon}
\newcommand{\ds}{\displaystyle}

\newcommand\calA{\mathcal{A}}
\newcommand\calD{\mathcal{D}}
\newcommand\calF{\mathcal{F}}
\newcommand\calJ{\mathcal{J}}
\newcommand\calS{\mathcal{S}}
\newcommand\calW{\mathcal{W}}
\newcommand\bbR{\mathbb{R}}
\newcommand\bbN{\mathbb{N}}
\newcommand\bbZ{\mathbb{Z}}

\newcommand{\ba}{\mbox{\boldmath$a$}}

\newcommand{\batilde}{\tilde{\ba}}
\newcommand{\bb}{\mbox{\boldmath$b$}}

\newcommand{\bbtilde}{\tilde{\bb}}
\newcommand{\bc}{\mbox{\boldmath$c$}}

\newcommand{\bctilde}{\tilde{\bc}}
\newcommand{\bd}{\mbox{\boldmath$d$}}
\newcommand{\bds}{{\mbox{\scriptsize \boldmath$d$}}}
\newcommand{\bdss}{{\mbox{\tiny \boldmath$d$}}}

\newcommand{\br}{\mbox{\boldmath$r$}}
\newcommand{\brtilde}{\tilde{\br}}

\newcommand{\bD}{\mbox{\boldmath$D$}}
\newcommand{\bDs}{{\mbox{\scriptsize \boldmath$D$}}}

\newcommand{\epsilontilde}{{\tilde{\epsilon}}}
\newcommand{\etabar}{{\bar{\eta}}}
\newcommand{\gammatilde}{{\tilde{\gamma}}}

\newcommand{\mutilde}{{\tilde{\mu}}}

\newcommand{\rhotilde}{{\tilde{\rho}}}
\newcommand{\sigmahat}{{\hat{\sigma}}}
\newcommand{\taubar}{{\bar{\tau}}}

\newcommand{\Phitilde}{\tilde{\Phi}}
\newcommand{\Psitilde}{\tilde{\Psi}}
\newcommand{\btilde}{{\tilde{b}}}
\newcommand{\atilde}{{\tilde{a}}}

\newcommand{\ctilde}{{\tilde{c}}}

\newcommand{\ktilde}{{\tilde{k}}}

\newcommand{\phat}{{\hat{p}}}
\newcommand{\ptilde}{{\tilde{p}}}
\newcommand{\rtilde}{{\tilde{r}}}

\newcommand{\Btilde}{{\tilde{B}}}

\newcommand{\Ctilde}{{\tilde{C}}}

\newcommand{\Dtilde}{{\tilde{D}}}
\newcommand{\Gbar}{{\bar{G}}}

\newcommand{\Tcheck}{{\check{T}}}
\newcommand{\That}{{\hat{T}}}
\newcommand{\What}{{\hat{W}}}

\newcommand{\Xbar}{{\bar{X}}}
\newcommand{\Xbbar}{{\Bar{\Bar{X}}}}
\newcommand{\Xhat}{{\hat{X}}}
\newcommand{\Ycheck}{{\check{Y}}}
\newcommand{\Yhat}{{\hat{Y}}}

\newcommand{\Zbar}{{\bar{Z}}}
\newcommand{\Zhat}{{\hat{Z}}}

\newcommand{\pp}[1]{\left({#1}\right)}
\renewcommand{\b}[1]{\left[{#1}\right]} 
\newcommand{\abs}[1]{\left\lvert#1\right\rvert}

\newcommand{\ind}[1]{\mathbbm{1}_{\{#1\}}}
\newcommand{\floor}[1]{\left\lfloor#1\right\rfloor}

\newcommand\cond{\, | \,}

\newcommand\comma{\, , \,}

\DeclareMathOperator{\bbE}{\mathbb{E}}
\DeclareMathOperator{\Var}{\mathrm{Var}}
\DeclareMathOperator{\bbP}{\mathbb{P}}

\newcommand\Deq{\stackrel{\mathrm{\calD}}{=}}
\newcommand\leqst{\stackrel{\mathrm{st}}{\leq}}
\newcommand\geqst{\stackrel{\mathrm{st}}{\geq}}
\newcommand{\etal}{{\it et al.}}

\newcommand{\dmax}{{d_{\mbox{\scriptsize max}}}}
\newcommand{\dmaxs}{{d_{\mbox{\tiny max}}}}

\newcommand{\toas}{\stackrel{\mathrm{a.s.}}{\to}}
\newcommand{\e}{{\mathrm e}}

\newcommand{\bbPD}{\bbP_{\bDs}}
\newcommand{\bbPDo}{\bbP_{\bDs(\omega_1)}}
\newcommand{\bbED}{\bbE_{\bDs}}
\newcommand{\bbEDo}{\bbE_{\bDs(\omega_1)}}

\newcommand{\dateenglish}{\renewcommand*{\today}{%
\number\day \ifcase\day \or
st\or nd\or rd\or th\or th\or th\or th\or th\or th\or th\or
th\or th\or th\or th\or th\or th\or th\or th\or th\or th\or
st\or nd\or rd\or th\or th\or th\or th\or th\or th\or th\or
st\fi\space \ifcase\month \or
January\or February\or March\or April\or May\or June\or July\or
August\or September\or October\or November\or December\fi \space\number\year}}

\newenvironment{remnn}{\addvspace{\medskipamount}\noindent\textbf{Remark.}}{\medbreak}
\newenvironment{remsnn}{\addvspace{\medskipamount}\noindent\textbf{Remarks.}}{\medbreak}

\numberwithin{equation}{section}

\newcommand{\hfigwidth}{6.5cm}
\newcommand{\figwidth}{11cm}

\begin{document}
\title{Threshold behaviour and final outcome of an\\epidemic on a random network with household\\ structure\footnote{First published in \emph{Advances in Applied Probability} 41:765--796 (2009); available at {\tt http://projecteuclid.org/euclid.aap/1253281063}.  Copyright (c) Applied Probability Trust 2009.}}

\author{
Frank Ball\thanks{School of Mathematical Sciences, University of Nottingham, University Park, Nottingham NG7 2RD, U.K. Email: {\tt frank.ball@nottingham.ac.uk}} \and
David Sirl\thanks{School of Mathematical Sciences, University of Nottingham, University Park, Nottingham NG7 2RD, U.K. Email: {\tt david.sirl@nottingham.ac.uk}} \and
Pieter Trapman\thanks{Julius Center for Health Sciences and Primary Care, University Medical Center Utrecht, Heidelberglaan 100, 3584 CX Utrecht, Netherlands and Stochastics Section, Faculty of Sciences, Vrije Universiteit Amsterdam, De Boelelaan 1081, 1081 HV Amsterdam, Netherlands. Current address: Department of Mathematics, Stockholm University, SE-10691 Stockholm, Sweden. Email: {\tt ptrapman@math.su.se}}
}

\date{13th May 2010}

\maketitle

\begin{abstract}
This paper considers a stochastic SIR (susceptible$\to$infective$\to$removed) epidemic model in which individuals may make infectious contacts in two ways, both within `households' (which for ease of exposition are assumed to have equal size) and along the edges of a random graph describing additional social contacts. Heuristically-motivated branching process approximations are described, which lead to a threshold parameter for the model and methods for calculating the probability of a major outbreak, given few initial infectives, and the expected proportion of the population who are ultimately infected by such a major outbreak. These approximate results are shown to be exact as the number of households tends to infinity by proving associated limit theorems. Moreover, simulation studies indicate that these asymptotic results provide good approximations for modestly-sized finite populations. The extension to unequal sized households is discussed briefly.
\end{abstract}


\section{Introduction}
\label{sec:Introduction}
Epidemic models which include some element of realistic population structure have been the subject of a considerable amount of recent study in recognition of the fact that the classical homogeneously-mixing models are quite unrealistic for all but the smallest of populations.

One approach to this has been to allow \emph{local} contacts of some kind, modelling contacts which occur on a regular basis in addition to maintaining the `well-mixed' \emph{global} contacts to model chance interactions with random members of the population. A common form for these local contacts to take arises by partitioning the population into \emph{households}, where these local contacts can occur only between individuals who are in the same household (see, for example, Becker and Dietz (1995) and Ball \etal\ (1997)). This can be extended to the overlapping groups model where the population may be partitioned in more than one way (for example, by household and by workplace), with local interactions taking place at (possibly) different rates within groups of the different partitions, see Ball and Neal (2002). Another mode of local interactions is described by the so-called great circle model (Ball \etal, 1997; Ball and Neal, 2002, 2003), where the population is spread around a circle and individuals have local contact with only their nearest neighbours. This model is closely related to `small-world' models (Watts and Strogatz, 1998), which have received considerable attention, particularly in the physics literature.

Another way of accounting for the inhomogeneous nature of interactions is by using random graphs to model social networks (see, for example, Andersson (1997, 1998, 1999), Newman (2002), Durrett (2006), Kenah and Robins (2007) and Britton \etal\ (2008)). Perhaps the most important aspect of these random graph models is that they incorporate a specified degree distribution, the degree of a node in the graph corresponding to the number of other members in the population an individual can possibly make infectious contact with. These models have been extended to also incorporate `casual contacts' by way of the classical homogeneous mixing effects, see Kiss \etal\ (2006) and Ball and Neal (2008).

In this paper we investigate a model for an SIR (susceptible$\to$infective$\to$removed) epidemic in a closed finite population, which draws together the main aspects of the generalisations of the standard homogeneously mixing model described above.   We consider a population grouped into households, with infectious contacts at a given per-pair rate, where individuals also make global contacts along the edges of a random graph over the whole population.  We use branching process approximations to derive (i) a threshold parameter, which determines whether a disease with just a few initial infectives can become established and infect a non-negligible proportion of the population (an event we call a \emph{major outbreak}); (ii) the probability that a major outbreak occurs; and (iii) the expected proportion of the population that is infected by a major outbreak.  These results are approximations that become exact in the limit as the size of the population becomes large in an appropriate way.

A feature of our model is that there is clustering present in the network of possible contacts, roughly meaning that there are significant numbers of triangles (and other short cycles) present in the network. This is an important aspect as the presence of triangles captures the phenomenon of people having mutual friends. The effect of such clustering in random networks in an epidemiological setting has been considered, in different models, by Trapman (2007) and Britton \etal\ (2008).

In the remainder of the paper we firstly describe, in Section~\ref{sec:Model}, the full detail of our model. Then in Section~\ref{sec:MainResults} we give the ideas behind the above-mentioned branching process approximations. In Section~\ref{sec:Calculations} we derive explicit formulae which allow us to calculate the quantities of interest for two important special cases, then give some brief numerical examples in Section~\ref{sec:Examples}, including demonstrating that our asymptotic results give good approximations for even moderately-sized finite populations. In Section~\ref{sec:Proofs} we rigorously establish the branching process approximations by proving related limit theorems as the population size tends to infinity. The paper concludes with a brief discussion in Section~\ref{sec:Discussion}. 

\section{Model}
\label{sec:Model}
We consider a closed population of $m$ households, each of $n$ individuals, and construct the network of possible global contacts using the `configuration model' (as in Durrett (2006, Chapter~3)) as follows. Firstly assign to each individual a number of \emph{half-edges}, these numbers being independent realisations of a random variable $D$ (the degree distribution) with $\bbP(D=k)=p_k,\,k=0,1,\ldots$. Conditional on the total number of half-edges being even we then pair these half-edges with each other uniformly at random, whence each such pair of half-edges forms an edge in the (random) graph describing the possible global contacts. We denote by $\mu_D$ and $\sigma^2_D$ the mean and variance of the distribution $D$ and assume that both of these quantities are finite. We also note for later reference that if we follow an edge from one vertex to another then the degree distribution of the second vertex is the \emph{size-biased} distribution $\Dtilde$, where $\bbP(\Dtilde=k)=kp_k/\mu_D,\, k=1,2,\ldots$. This is because in the construction of the graph the half-edges are paired uniformly at random, so it is $k$ times more likely that following an edge leads one to a vertex of degree $k$ than to a vertex of degree $1$. By the degree of an individual we mean the number of individuals adjacent to it in the network of global contacts, not counting those in its own household.

Note that there may be some imperfections in the graph, in the form of parallel edges and self-loops. However, our assumption that $\sigma^2_D<\infty$ ensures that as $m\to\infty$, the number of these imperfections in the network of global contacts converges in distribution to a Poisson random variable whose mean is a function of $(\mu_D,\sigma^2_D)$ (Durrett, 2006, Theorem~3.1.2). By treating the households as macro-individuals, with degree distribution given by the sum of $n$ independent copies of $D$, it follows that the numbers of parallel edges between households and household self-loops also converge in distribution to Poisson random variables as $m\to\infty$. Thus the probability that these imperfections are absent in the graph is bounded away from zero as $m \to \infty$, and consequently (cf.~Janson (2009)) our asymptotic results continue to hold if the graph is conditioned on having no such imperfections.

When an infective individual makes infectious contact with a susceptible individual, the susceptible becomes infective and remains so for a random period of time distributed according to a non-negative real-valued random variable $I$, which we specify by its Laplace transform $\phi(\theta)=\bbE [\e^{-\theta I}],\, \theta\geq0$, and call the infectious period. An infective individual makes infectious contact with each other member of his/her household at the points of a Poisson process with rate $\lambda_L$ and similarly with each individual he is adjacent to in the network of global contacts at rate $\lambda_G$. To be emphatic, both $\lambda_L$ and $\lambda_G$ are \emph{per-pair} rates, so an infectious individual of degree $k$ makes infectious contacts at overall rate $\lambda_L(n-1)+\lambda_Gk$. As usual, all Poisson processes and infectious periods are assumed to be mutually independent.

For ease of presentation, we assume that an epidemic is initiated by a single infective individual within the population, either a given specific individual or an individual chosen uniformly at random from the population. Our assumption that all households are of the same size is also made for ease of presentation although, as indicated in Section~\ref{sec:Discussion}, our results generalise easily to incorporate unequal household sizes.

\section{Heuristics and description of main results}
\label{sec:MainResults}

We now give informal descriptions of the branching process approximations we use, firstly to approximate the early stages of an outbreak, leading to a threshold parameter and a method of calculating the probability of a major outbreak and, secondly, to approximate the expected relative final size of (i.e.\ the proportion of the population infected by) a major outbreak.  These approximations become exact in the limit as the number of households $m \to \infty$, with the household size $n$ held fixed.

\subsection{Forward processes}
\label{sec:MainFwd}
The branching process we use to analyse the early stages of the epidemic approximates the number of households which become infected in the course of the epidemic. Because we are interested only in the final outcome of the epidemic and not its precise time evolution we can think of the epidemic as evolving in the following way (see, for example, Pellis \etal~(2008)). We first consider the epidemic spreading only within the household containing the initial infective (the local epidemic that it initiates) and then consider the number of individuals infected via global infectious contacts made by those infected by the local epidemic. Because of the way the network is constructed, in the early stages of the epidemic it is highly likely that these globally contacted individuals are all in distinct households (this being critical for the branching process approximation). We then consider each newly infected household in the same manner: local epidemic followed by global infections. Again in the early stages it is highly likely that those infected by such global infectious contacts are in distinct households and furthermore that they are in previously uninfected households. We can view this as a branching process if we consider the households infected by a local epidemic initiated by a single infective within a typical household to be the children (offspring) of that household.

Note that the offspring distribution of the above branching process is different for the initial (i.e.\ zeroth) generation than for subsequent generations, since in subsequent generations the initial infective in a household has been infected by one of its global neighbours, so the number of uninfected global neighbours of this individual is equal in distribution to $\Dtilde-1$, whilst in the zeroth generation the initial infective is the initial infective in the whole population, and the degree distribution of this individual is either distributed as $D$ or is a fixed constant, according as the initial infective is chosen (uniformly) at random or a specific individual in the population is chosen to be the initial infective. We therefore define the random variable $C$ to be the number of global neighbours infected by members of the initial infective's household and $\Ctilde$ to be the number infected by the household of a single infective that was infected by a global neighbour. Our branching process approximation is then defined by it having a single ancestor (since the epidemic starts with one initial infective) and offspring distribution $C$ in the initial generation and $\Ctilde$ in subsequent generations.  Throughout the paper, we denote a branching process of this type by ${\rm BP}(1,C,\Ctilde)$, or by ${\rm BP}(1,\bc,\bctilde)$, where $\bc=(c_0,c_1,\ldots)$ and $\bctilde=(\ctilde_0,\ctilde_1,\ldots)$ are the mass functions of $C$ and $\Ctilde$, respectively.

The above branching process approximation of the epidemic is made fully rigorous in Section~\ref{sec:threshold}, where it is shown that, as $m \to \infty$, the total number of households infected by the epidemic converges in distribution to the total progeny of the branching process (see Theorem~\ref{thm:threshold}).  Thus, whether or not the epidemic can `take off' and lead to a major outbreak is determined by whether or not the branching process is supercritical (i.e.~whether or not $R_*=\bbE [\Ctilde]>1$).  Further, by standard branching process theory, the probability of such a major outbreak is given by $1-f_C(\sigma)$, where $\sigma$ is the smallest solution of $f_\Ctilde(s)=s$ in $[0,1]$, and $f_C(s)=\bbE [s^C]$ and $f_\Ctilde(s)=\bbE [s^\Ctilde]$ (for $s\in[0,1]$) denote the probability generating functions (PGFs) of $C$ and $\Ctilde$, respectively. (Here and henceforth we denote by $f_X(\cdot)$ the PGF of the random variable $X$.)  Calculation of $R_*$, $f_C(s)$ and $f_\Ctilde(s)$ is considered in Section~\ref{sec:FwdCalculations}.

\subsection{Backward processes}
\label{sec:MainBwd}
We now consider the expected final size of a major outbreak. Again our analysis is of the $m\to\infty$ limiting epidemic process, for which we find the probability that a given individual is infected in the event of a major outbreak. By an exchangeability argument this probability is equal to the asymptotic mean proportion of the population (individuals, not households) that are ultimately infected by a major outbreak. This quantity serves as our approximation of the expected proportion infected in a major outbreak in a finite population. We determine the probability that a given individual is infected by considering its \emph{susceptibility set} (cf.~Ball and Lyne (2001) and Ball and Neal (2002)).

The idea behind susceptibility sets is that for each individual in the population we can, by sampling from the infectious period distribution and then the relevant Poisson processes, make a (random) list of other individuals it would infect were it to be infected itself. We then construct a digraph (directed graph) based on these lists, in which the vertices represent individuals in the population and we put a directed arc from $i$ to $j$ when, were $i$ to become infected, it would make infectious contact with $j$, i.e.~if $j$ is in $i$'s list. The susceptibility set of individual $i$ consists of those individuals from which there exists a path to $i$ in the digraph (including $i$ itself). Note that an individual will become infected by an epidemic if and only if the initial infective is in its susceptibility set.  We also need the concept of a \emph{local} susceptibility set, constructed in the same way but considering only local (within-household) infectious contacts.

We approximate the size of the susceptibility set of an individual chosen uniformly at random from the population by the total progeny of an appropriate branching process. To construct this branching process we break up the susceptibility set into `generations' in much the same way as we look at the spread of infection in the early stages of the epidemic. Starting with an individual $i$, consider those individuals $j$, not in $i$'s household, who are in $i$'s susceptibility set by virtue of an arc leading from $j$ to an individual in $i$'s local susceptibility set. These individuals are all in different households with high probability as $m\to\infty$ and the households they are in comprise the first `generation' of the susceptibility set. Repeating this process for each of these individuals $j$ (i.e.~looking at the individuals who make infectious global contact with a member of $j$'s local susceptibility set) gives the second `generation'; and by continuing this process we can construct the whole of $i$'s susceptibility set.
Because each individual $j$ that joins the susceptibility set by virtue of a global contact is in a household not previously associated with the susceptibility set with high probability, the number of households in each generation is approximated well by the branching process ${\rm BP}(1,B,\Btilde)$, where $B$ and $\Btilde$ denote the offspring random variables for the initial and subsequent generations, which again are typically different.

We show in Section~\ref{sec:majOBInfProp} that, as $m\to\infty$, the conditional probability that a typical initial susceptible ($i$ say) is infected, given that a major outbreak occurs, is given by the probability that the branching process ${\rm BP}(1,B,\Btilde)$ avoids extinction (see Theorem~\ref{thm:majOBInfProb}). An intuitive explanation of this result is as follows.  As $m\to\infty$, (i) the number of households in $i$'s susceptibility set converges in distribution to the total progeny of ${\rm BP}(1,B,\Btilde)$; and (ii) a major outbreak necessarily infects at least $\log m$ households (cf.~Lemma~\ref{lem:limProbYt}).  Thus, as $m \to \infty$, the probability that $i$'s susceptibility set intersects one of these $\log m$ households is 0 if ${\rm BP}(1,B,\Btilde)$ goes extinct and 1 otherwise.  The latter result follows because if ${\rm BP}(1,B,\Btilde)$ does not go extinct then the size of $i$'s susceptibility set is of exact order $m$ as $m \to \infty$.

The above claim and standard branching process theory imply that the expected relative final size of a major outbreak in a large finite population is approximately $1-f_B(\xi)$, where $\xi$ is the smallest solution of $f_\Btilde(s)=s$ in $[0,1]$.  Calculation of $f_B(s)$ and $f_\Btilde(s)$ is considered in Section~\ref{sec:BwdCalculations}.

\section{Calculations}
\label{sec:Calculations}

\subsection{Forward process}
\label{sec:FwdCalculations}

Consider first the threshold parameter $R_* = \bbE [\Ctilde]$. Label the individuals in a household $0,1,\ldots,n-1$, with individual $0$ the initial infective, and define $\chi_i$ to be the indicator of the event that individual $i$ is infected in the local (i.e.~single-household) epidemic and $C_i$ to be the number of global neighbours with which $i$ makes infectious contact, if $i$ were to become infected.  Then
\begin{equation}
\Ctilde=C_0 + \sum_{i=1}^{n-1}\chi_i C_i
\label{eq:C}
\end{equation}
and it follows, since $C_1$ and $\chi_1$ are independent and $(C_1,\chi_1),(C_2,\chi_2),\ldots,(C_{n-1},\chi_{n-1})$ are identically distributed, that
\begin{equation}
R_* = \bbE [C_0] + \bbE [T] \bbE [C_1],
\label{eq:Rstar2}
\end{equation}
where $T$ is the final size of the within-household epidemic (not counting the initial infective). Denote by $I_i$ and $K_i$ the infectious period and number of global neighbours, not including its infector, of individual $i$ (this only affects the initial infective within the household). Now, since infectious contacts between different pairs of individuals are independent, $C_i \cond K_i , I_i \sim \mbox{Bin}(K_i, 1-\e^{-\lambda_G I_i})$. Thus $\bbE [C_i \cond K_i , I_i] = K_i(1-\e^{-\lambda_G I_i})$, whence, by the independence of $K_i$ and $I_i$,
\begin{equation}
\bbE [C_i] = \bbE [K_i] (1-\phi(\lambda_G)).
\label{eq:ECi}
\end{equation}
Now, for $i=1,2,\ldots,n-1$, $K_i$ has the same distribution as $D$, the prescribed degree distribution, so $\bbE [K_i]=\mu_D$ for such $i$. However, for the reasons noted in the first paragraph of Section~\ref{sec:Model}, since the initial infective in the household was infected by a global infection its degree has the size-biased distribution $\Dtilde$, and because one of these neighbours (the one that infected it) has already been infected, $K_0$ has the same distribution as $\Dtilde-1$, so $\bbE [K_0] = \bbE [\Dtilde]-1$. It follows from the definition of $\Dtilde$ that $\bbE [\Dtilde] = \bbE [D] + \Var D/\bbE [D]$. Substituting these into~\eqref{eq:ECi} and then~\eqref{eq:Rstar2}, and letting $\mu_T=\bbE [T]$, yields
\begin{equation}
R_* = \pp{\mu_D \: (\mu_T + 1) + \frac{\sigma^2_D}{\mu_D} -1 } (1-\phi(\lambda_G)).
\label{eq:RstarFinal}
\end{equation}
The mean $\mu_T$ may be evaluated (typically numerically) by using equations~(2.25) and~(2.26) of Ball (1986), thus enabling $R_*$ to be calculated.

Calculation of the PGFs $f_C(s)$ and $f_\Ctilde(s)$ is more difficult because the number of global infections caused by a particular individual is dependent on that individual's infectious period, which also influences whether or not other individuals in the household become infective and thus the number of global contacts they might make. It is possible to use the notion of `final state random variables' introduced by Ball and O'Neill (1999) to find $f_C$ and $f_{\Ctilde}$, but it is not straightforward, so we do not present it here. This methodology will be discussed in a forthcoming paper concentrating on the more applied aspects of our model. However, there are two special cases where the above dependencies do not exist and the analysis is much simpler. These are when the infectious period is fixed (i.e.~almost surely equal to a given constant) and when the infectious period can be only zero or infinity.

Trapman (2007) describes (using results of Kuulasmaa (1982)) how these special cases lead to bounds on quantities of interest for a very general class of epidemic models. Trapman's arguments hold for any epidemic model where there is only one `kind' of infectious contact rather than the two (local and global) that we are concerned with, but the methods can be easily adapted. In addition, a fixed infectious period is often a reasonable assumption to make in practice and it is commonly used because it leads to simplifications of the kind shown shortly (see, for example, Britton \etal\ (2007) and Britton \etal\ (2008)). We therefore proceed to calculate the PGFs $f_C$ and $f_{\Ctilde}$ in these two special cases as they can be used to calculate the above-mentioned bounds and they also may give insight into the importance of and interplay between the parameters of our model. The role of the infectious period distribution will be discussed in the above-mentioned applied paper.

\subsubsection{Zero or infinite infectious period.}
Suppose that $\bbP(I=\infty)=1-\bbP(I=0)=p$ for some $p\in[0,1]$. For the moment we ignore the differences between the initial and subsequent generations and denote the generic offspring random variable by unadorned $C$. Here we have
\begin{equation*}
C=\begin{cases}
0, & \mbox{with probability $1-p$}, \\
C_0 + \sum_{i=1}^{n-1} C_i, & \mbox{with probability $p$},
\end{cases}
\end{equation*}
where $C_i$ is the number of global neighbours infected by an infectious individual $i$. Thus $C_0 \Deq K_0$ (where $\Deq$ denotes equality in distribution) and $C_1,C_2,\ldots,C_{n-1}$ are independent and identically distributed with
\begin{equation*}
C_i=\begin{cases}
0, & \mbox{with probability $1-p$}, \\
K_i, & \mbox{with probability $p$}.
\end{cases}
\end{equation*}
Also note that the number, $N$ say, of the $n-1$ $C_i$'s which take the value $K_i$ (i.e.~the number of initially susceptible individuals in the household with $I=\infty$) is binomially distributed, with parameters $n-1$ and $p$. We therefore have
\begin{eqnarray*}
f_C(s)=\bbE [s^C] & = & (1-p)s^0 + p\bbE [s^{C_0 + \sum_{i=1}^{n-1} K_i}] \\
 & = & 1-p + p\bbE [s^{C_0}] \bbE [s^{\sum_{i=1}^N K_i}] \\
 & = & 1-p + p f_{K_0}(s) f_N(f_D(s)) \\
 & = & 1-p + p f_{K_0}(s) (1-p+pf_D(s))^{n-1},
\end{eqnarray*}
where $K_0$ is $D$ or $d$ in the initial generation and $\Dtilde-1$ in subsequent generations (in which case the PGF is $f_{\Ctilde}$ rather than $f_C$).

\subsubsection{Fixed infectious period.}
Now suppose that $\bbP(I=c)=1$ for some $c>0$. Again we temporarily ignore the differences between the initial and subsequent generations, label the individuals $0,1,\ldots,n-1$ and denote by $C_i$ the number of global neighbours infected by an infectious individual $i$. Then, letting $T$ denote the final size of the within-household epidemic, we have $C=C_0+\sum_{i=1}^T C_i$ and, conditional on the final size, $C_1,C_2,\ldots,C_T$ are mutually independent. Now $C_i \cond K_i \sim \mbox{Bin}(K_i,1-\e^{-c\lambda_G})$, so $f_{C_i}(s)=f_{K_i}(1-p_G+sp_G)$, where $p_G=1-\e^{-c\lambda_G}$. Thus, by the usual formula for the PGF of a random sum,
\begin{equation}
\label{eq:fcsfixed}
f_C(s) = f_{C_0}(s) f_T(f_{C_1}(s)) = f_{K_0}(1-p_G+sp_G) f_T(f_D(1-p_G+sp_G)),
\end{equation}
where again $K_0$ is $D$ or $d$ in the initial generation and $\Dtilde-1$ in subsequent generations. The PGF $f_T$ is easily calculated using Theorem~2.6 of Ball (1986).

\subsection{Backward process}
\label{sec:BwdCalculations}
Now consider the branching process approximation of the growth (as described in Section~\ref{sec:MainBwd}) of the susceptibility set of an individual, $i_*$ say, chosen uniformly at random from the population. The offspring distribution of this process has the same distribution as the number of individuals that make global contact with the local susceptibility set of a single individual, say individual $i$. Again we have a distinction between the initial and subsequent generations but we ignore this for now and denote the random variable of interest by $B$. Firstly we write
\begin{equation*}
B=B_0 + \sum_{j=1}^M B_j,
\end{equation*}
where $B_j$ is the number of contacts made with individual $j$ (again labelling the individuals within the household $0,1,\ldots,n-1$, with $0$ corresponding to the primary individual $i$) and $M$ is the size of $i$'s local susceptibility set, not counting $i$ itself. (If $M=0$ then $i$'s local susceptibility set consists of only $i$ itself and the sum is empty.) Now $B_j \cond K_j \sim\mbox{Bin}(K_j,p_G)$, where $K_j$ is the number of global neighbours of $j$ excluding, in the case of the initial individual, the individual it made contact with in order to join the susceptibility set and $p_G=1-\phi(\lambda_G)$ is the probability that an infective individual makes infectious contact with a given global neighbour. We do not need to condition on the infectious period of individual $j$ because the contacts we are considering come from other individuals; the independence of the infectious periods of these individuals implies that they make contacts with $j$ independently of each other. For a similar reason, $B_0,B_1,\ldots,B_M$ are independent.  Arguing as in the the derivation of \eqref{eq:fcsfixed} yields that
\begin{equation}
\label{eq:fbs}
f_B(s)=f_{K_0}(1-p_G+sp_G) f_M(f_{D}(1-p_G+sp_G)),
\end{equation}
where now $K_0$ is $D$ in the initial generation (because of how $i_*$ was chosen) and $\Dtilde-1$ in subsequent generations.

In order to determine $f_M$ we use equation~(3.5) of Ball and Neal (2002), which gives a triangular system of linear equations whose solution is the mass function of $M$, from which one can easily calculate the PGF.  Note that \eqref{eq:fbs} holds for any choice of infectious period distribution.  It is easily verified that in the fixed infectious period case $T\Deq M$, so $f_{\Btilde}(s)=f_{\Ctilde}(s)$ and, if the initial infective is chosen uniformly at random from the population, $f_B(s)=f_C(s)$; and in the zero or infinite infectious period case $M\sim \mbox{Bin}(n-1, p)$, where $p=\bbP(I=\infty)$, whence $f_M(s)=(1-p+ps)^{n-1}$.

\section{Numerical results}
\label{sec:Examples}
\begin{figure}[h!]
\begin{center}
\psfrag{lamG}[][]{\small $\lambda_G$}
\psfrag{lamL}[][]{\small $\lambda_L(n-1)$}
\includegraphics[width=\figwidth]{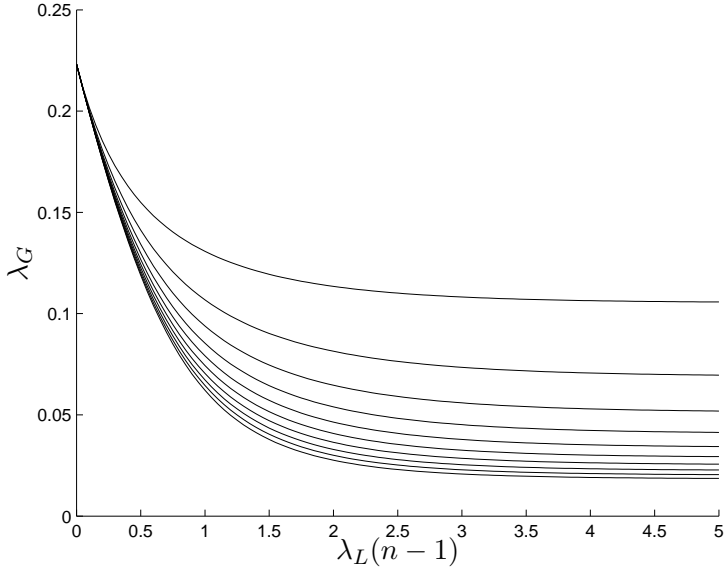}
\caption{Critical values of $\lambda_G$ and $\lambda_L(n-1)$ above which the epidemic is supercritical, for $n=2,3,\ldots,10$ (top to bottom in the plot). Other parameters are $I\equiv1$ and $D\sim\mbox{Poi}(5)$ (i.e.~Poisson with mean 5).}
\label{fig:critLamGplot}
\end{center}
\end{figure}

We now explore, numerically, some of the features of our model and investigate how they depend on some of its parameters. As a way of examining how the household size $n$ affects the model, Figure~\ref{fig:critLamGplot} shows the critical values of the per-pair global contact rate $\lambda_G$ and the per-individual local contact rate $\lambda_L(n-1)$, above which the epidemic is supercritical, for several household sizes, with the degree distribution and infectious period distribution fixed. Note that the expected total rate of global contacts per individual remains constant over these plots since $D$ is held fixed.
Note also that if $\lambda_L=0$ then $n$ is immaterial, as is $\lambda_L$ when $n=1$. In these situations there is no local contact, so we recover the standard network model and the critical value of $\lambda_G$ is at the point the plotted lines converge to as $\lambda_L \to0$. The plot reflects the fact that, even as the per-individual total contact rate remains constant, increasing the household size spreads the potential infectious contacts over a larger number of neighbours, thus avoiding repeated contacts with the same individual and increasing the spread of the disease. We also observe that, fixing $D$ and letting $\lambda_L\to\infty$, the critical value of $\lambda_G$ tends to that for the standard network model with the same infectious period distribution and degree distribution $\sum_{i=1}^n D_i$, where the $D_i$ are independent copies of $D$. This is because, in this limit, once an individual is infected the whole household that it is in necessarily becomes infected, and is easily verified using~\eqref{eq:RstarFinal}.

\begin{figure}[h!]
\begin{center}
\psfrag{probMajOB}[][]{\small Major outbreak probability}
\psfrag{muD}[][]{\small $\mu_{{\!}_D}$}
\psfrag{Poisson}{\scriptsize Poisson}
\psfrag{GeometricLong}{\scriptsize Geometric}
\psfrag{Constant}{\scriptsize Constant}
\psfrag{Hlabel}{\scriptsize Power law}
\psfrag{HClabel}{\scriptsize Power cutoff}
\includegraphics[width=\figwidth]{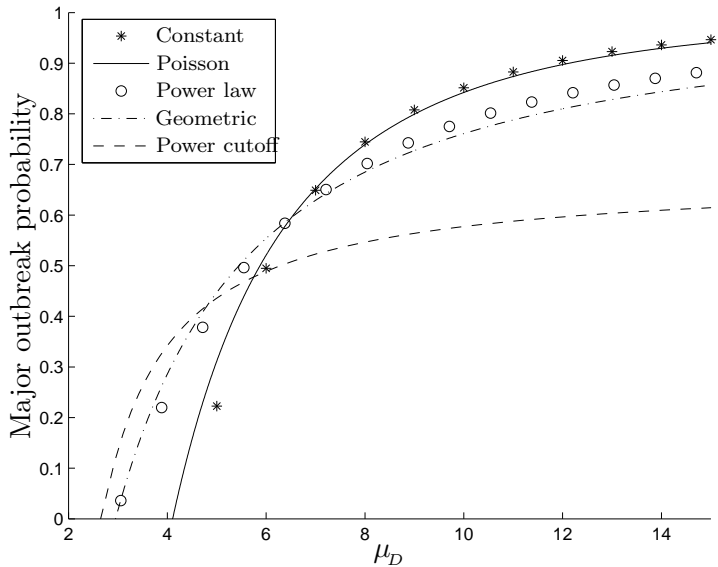}
\caption{The probability of a major outbreak versus $\mu_D$ for different classes of degree distribution $D$. The distribution labelled `Power law' is $\mbox{Pow}(k_*,7/2)$, for $k_*=5,6,\ldots,18$ and the distribution labelled `Power cutoff' is $\mbox{PowC}(\kappa,3/2)$, for $\kappa\in[10,485]$ (smaller values of $k_*$ or $\kappa$ yield subcritical epidemics). The other parameters of the model are $n=3$, $I\equiv1$, $\lambda_L=1$ and $\lambda_G=1/10$.}
\label{fig:MajOBProbVsD}
\end{center}
\end{figure}

Perhaps the most interesting aspect of this model to explore is the dependence of its behaviour on the distribution of $D$, the number of global neighbours of a typical individual. Considerable research, conjecture and discussion has gone into trying to determine distributions which capture the features of many real life contact networks---Section~III.C of Newman (2003) has an extensive list of references. In Figure~\ref{fig:MajOBProbVsD} we investigate the probability of a major outbreak in our epidemic model for various distributions $D$ with different properties, in particular different tail behaviours. We use the standard Poisson and geometric (with support including 0) distributions, as well as an almost surely constant degree and two variants of heavy-tailed distributions. The first has mass function
\begin{equation*}
p_k \propto \begin{cases}
k_*^{-a}, & \mbox{for $k=1,2,\ldots,k_*$}, \\
k^{-a}, & \mbox{for $k=k_*+1, k_*+2,\ldots$,}
\end{cases}
\end{equation*}
and the second, with mass function $p_k \propto k^{-a}\e^{-k/\kappa}$ ($k=1,2,\ldots$), is a power law with exponential cut-off which has gained much attention in recent physics literature. We denote these distributions by $\mbox{Pow}(k_*,a)$ and $\mbox{PowC}(\kappa,a)$, respectively.
The behaviour of these plots for relatively small values of $\mu_D$ (where the model is close to critical) is largely determined by the probability of $D$ taking very large values, i.e.~the tail of the distribution, as this dictates what opportunity the disease might have to really `take hold'; however when $\mu_D$ is large the behaviour of $D$ at small values is more important, as the epidemic can usually move quite freely and this determines the chance that it might be contained by the network structure.

\begin{figure}
\begin{center}
\psfrag{titlePMOPoi}[][]{\tiny (a) $D\sim\mbox{Poi}(8)$}
\psfrag{titleRFSPoi}[][]{\tiny (b) $D\sim\mbox{Poi}(8)$}
\psfrag{titlePMOHeavy}[][]{\tiny (c) $D\sim\mbox{Pow}(10,7/2)$}
\psfrag{titleRFSHeavy}[][]{\tiny (d) $D\sim\mbox{Pow}(10,7/2)$}
\psfrag{nohh}[][]{\tiny Number of households}
\psfrag{PMO}[][]{\tiny Major outbreak probability}
\psfrag{RFS}[][]{\tiny Relative final size}
\includegraphics[width=\hfigwidth]{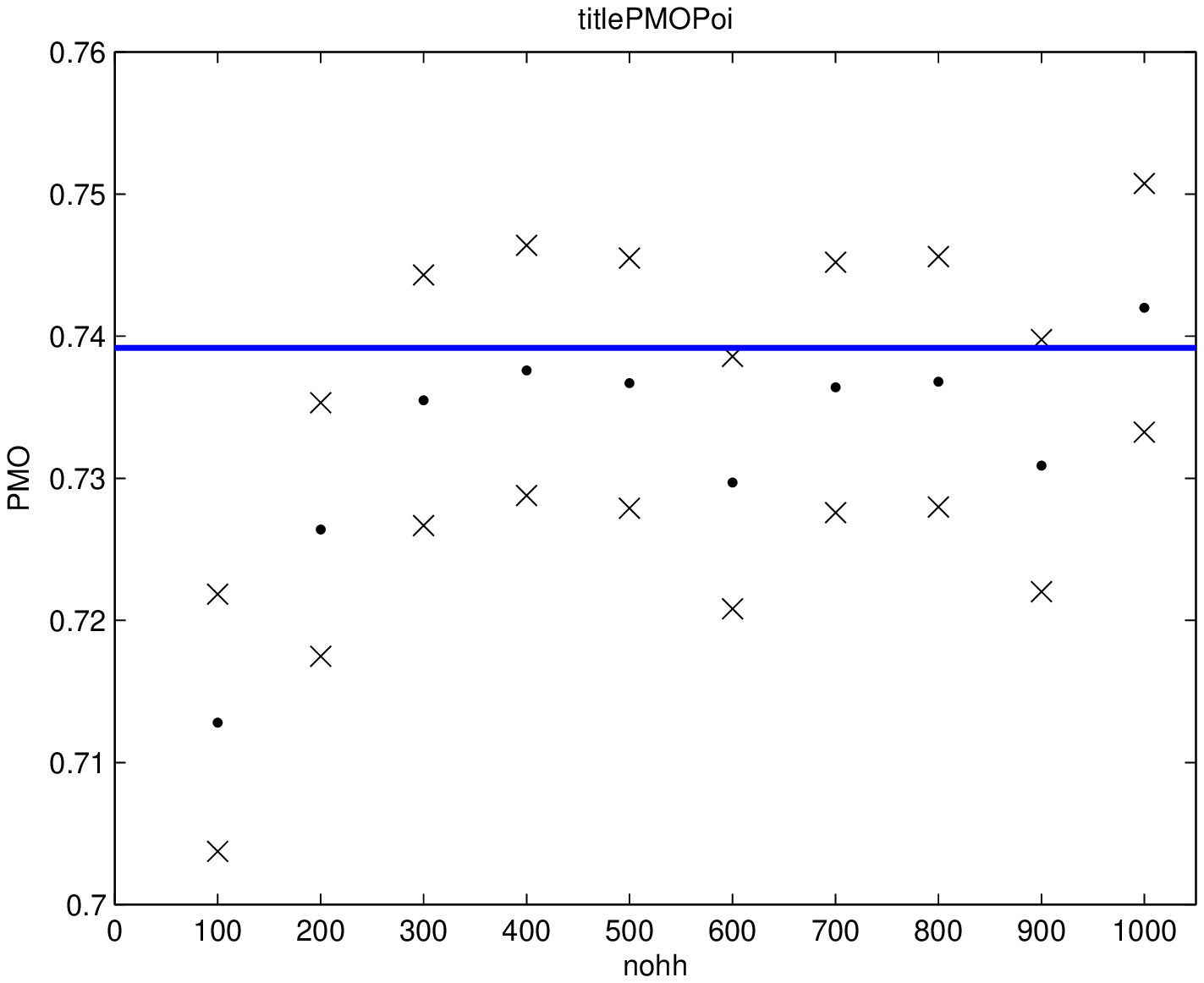}
\includegraphics[width=\hfigwidth]{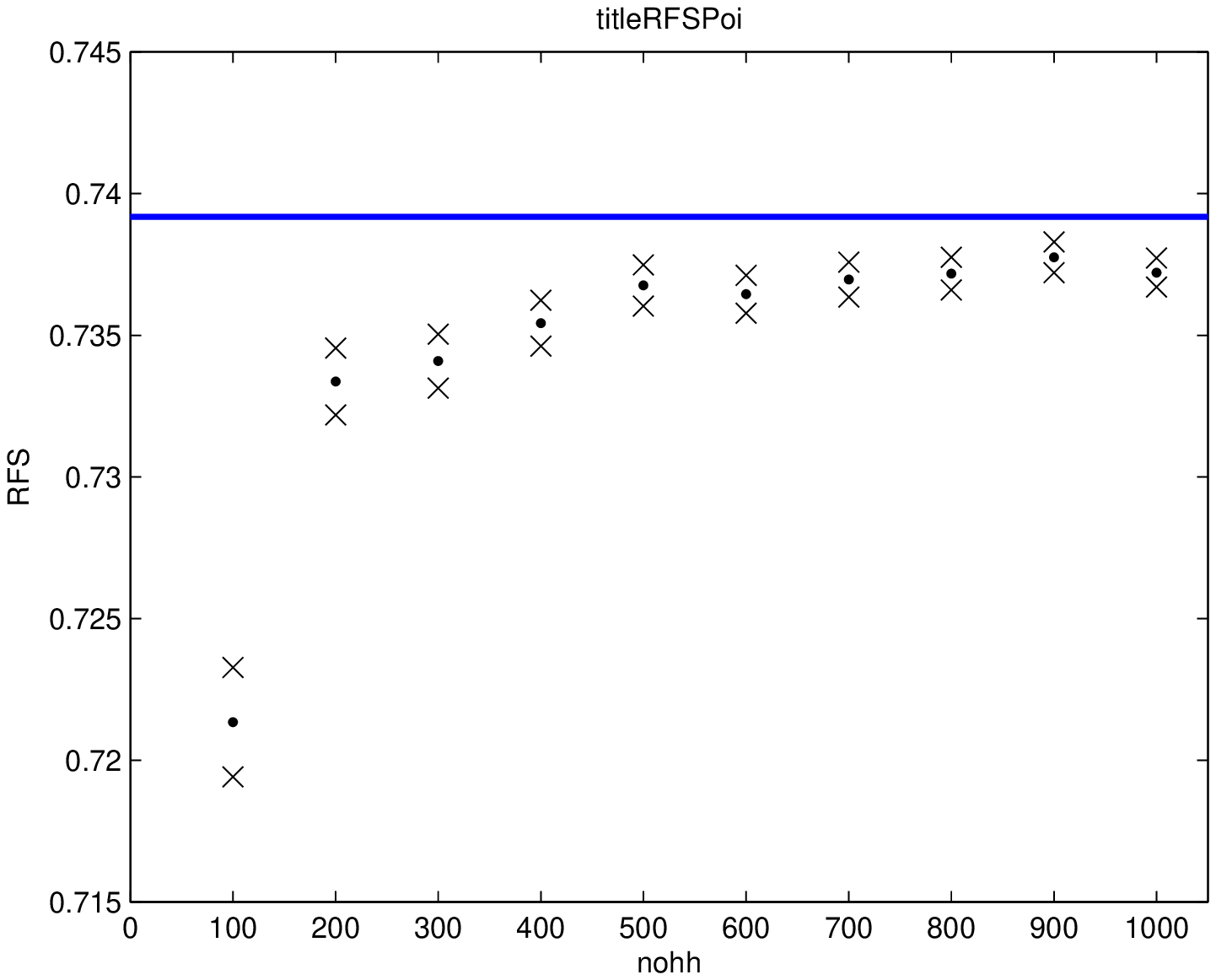} \\
\includegraphics[width=\hfigwidth]{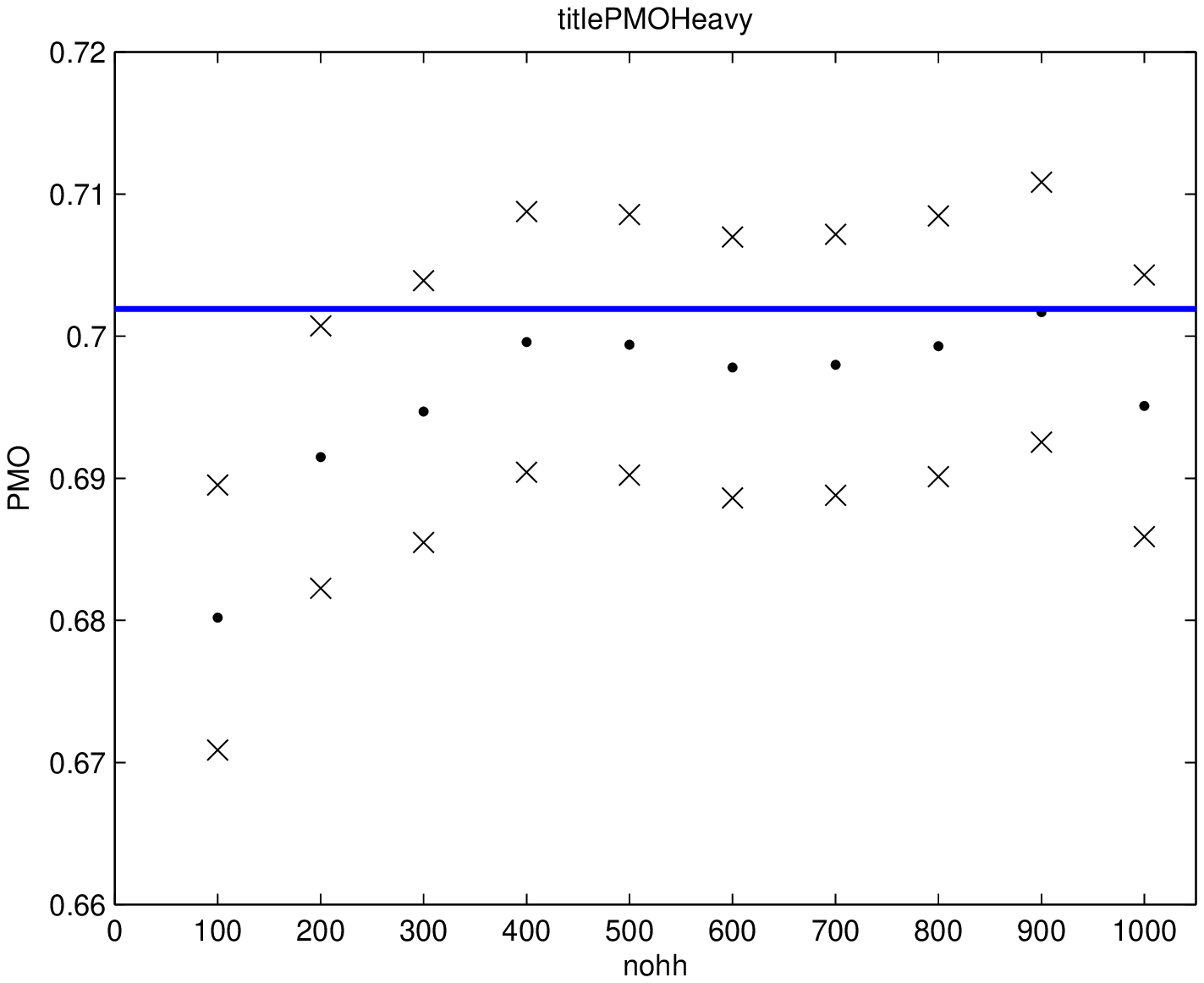}
\includegraphics[width=\hfigwidth]{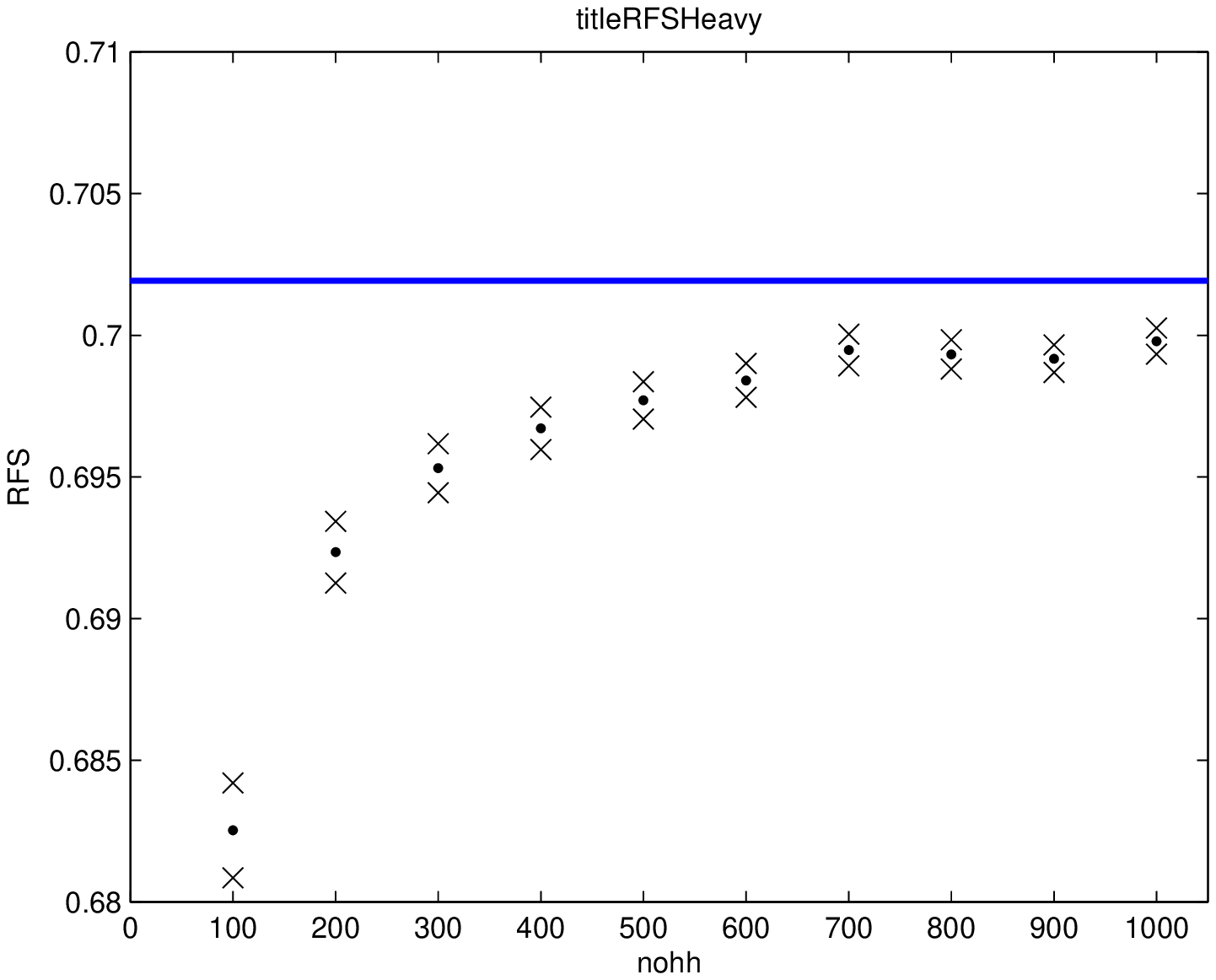}
\end{center}
\caption{Comparison of simulation estimates of major outbreak probability and expected relative final size for finite populations with asymptotic results. The Poisson degree distribution (plots (a) and (b)) has $\mu_D=\sigma^2_D=8$ and the power law distribution (plots (c) and (d)) has $\mu_D\approx8.04$ and $\sigma^2_D\approx96$. Other parameters are $n=3$, $I\equiv1$, $\lambda_L=1$ and $\lambda_G=1/10$.}
\label{fig:cgce}
\end{figure}

We also briefly investigate whether our asymptotic methods give reasonable approximations to the quantities of interest in finite populations. We estimate the probability and expected relative final size of a major outbreak in finite populations from simulations and compare these to the results we get from our asymptotic analysis. Each simulation consists of generating a random network and running \emph{one} epidemic on it. Figure~\ref{fig:cgce} shows estimates of these quantities of interest for increasing numbers $m$ of households together with the theoretical ($m=\infty$) values for two choices of degree distribution. The estimates of the major outbreak probability are based on 10,000 simulations for each parameter combination and those that result in a major outbreak are then used to estimate the expected relative final size. We have plotted point estimates of the quantities of interest, together with error bounds based on $\pm2$ standard errors (SE) of the estimator. For the probability of a major outbreak, estimated as $\phat$, $\mbox{SE}= [\phat(1-\phat)/n_0]^{1/2}$, where $n_0=10,000$ is the number of simulations. For the relative final size, $\mbox{SE}= \sigmahat n_1^{-1/2}$, where $\sigmahat^2$ is the sample variance of the relative final sizes and $n_1$ is the number of simulations that resulted in a major outbreak.

Note that in small finite populations the determination of a cutoff for whether a particular final size constitutes a major outbreak is practically impossible; only once the population size is sufficiently large (for $m$ larger than about 100 in our simulations) does the distinction become clear. In our calculations we have used a cutoff of 0.15 of the population size, this being determined by inspecting histograms of the relative final size of the simulations. Also note that the vertical scale of plots (a) and (c) is different from that of plots (b) and (d). Figure~\ref{fig:cgce} shows that our asymptotic results give good approximations for these quantities of interest for populations of only a few hundred households. Though the asymptotic values of both the major outbreak probability and the expected relative final size seem to consistently overestimate these values for the finite populations (as one would expect since the approximating branching process treats each global infection as an infection of a previously uninfected household, thus overestimating disease spread), even for populations of only 100 households the relative error is much less than 5\%. It also seems that having a heavy-tailed degree distribution may make the convergence to the asymptotic value a little slower (compare plots (b) and (d) at around 200--500 households), but the effect seems to be only very slight. Another interesting observation is that the relative final size seems to be appreciably more efficiently estimated by our simulation methods than the probability of a major outbreak. This is owing (at least in part) to the fact that from each simulation we simply observe the occurrence or otherwise of a major outbreak---one observation of the forward process---whereas when a major outbreak does occur, the proportion infected has information about the susceptibility set of every initial susceptible in the population---many (highly correlated) observations of the backward process.

\section{Proofs}
\label{sec:Proofs}
\subsection{Overview}
In this section we provide a fully rigorous justification of the results discussed in Section~\ref{sec:MainResults} concerning the threshold behaviour of the epidemic model and its final outcome in the event of a major outbreak. This subsection gives a brief outline of our methods of proof. The starting point is a sequence $\bD=(D_1,D_2,\ldots)$ of independent copies of $D$. For $m=1,2,\ldots$, $(D_1,D_2,\ldots,D_{mn})$ is used to give the degrees of the $mn$ individuals in a population of $m$ households. We then define a realisation of the epidemic, $E^{(m)}$ say, viewed on a generation basis, and a realisation of an approximating branching process, say $Y^{(m)}=(Y^{(m)}_k,\, k=0,1,\ldots)$ (see Section~\ref{sec:construction}). In $E^{(m)}$ the network is formed, i.e.~the half-edges are paired up, as the epidemic progresses. The branching process $Y^{(m)}$ is similar to the branching process, $Y$ say, described in Section~\ref{sec:MainFwd}, except the empirical distribution of the degrees $D_1,D_2,\ldots,D_{mn}$ is used in place of the degree distribution $D$. The epidemic $E^{(m)}$ and approximating branching process $Y^{(m)}$ are coupled so that they coincide until a random number, $\tau^{(m)}+1$, of households have been infected in $E^{(m)}$. It is shown that $\bbP(\tau^{(m)}>k) \to 1$ as $m\to\infty$ for all $k\in\bbZ_+$, so $\Zhat^{(m)}$, the number of households infected in $E^{(m)}$, and $\Yhat^{(m)}$, the total progeny of $Y^{(m)}$, have the same limiting distribution as $m\to\infty$. (We use $\bbZ_+$ to denote the positive integers including 0 and $\bbN$ to denote the strictly positive integers.) Now, $Y^{(m)}$ converges in distribution to $Y$ as $m\to\infty$, so $\Zhat^{(m)}$ is asymptotically distributed as $\Yhat$, the total progeny of $Y$ (see Theorem~\ref{thm:threshold}), thus providing a formal justification of the threshold behaviour described in Section~\ref{sec:MainFwd}.

Suppose now that $R_*>1$, so that major outbreaks are possible. Let $t_m=\floor{2\log\log m/\log R_*}$, where, for $x\in\bbR$, $\floor{x}$ denotes the greatest integer $\leq x$. We show (cf.~Lemma~\ref{lem:limProbYtildeThat}) that there exists $\beta>1$ such that $\lim_{m\to\infty} \bbP(\log m < Z^{(m)}_{t_m} < (\log m)^\beta) = \bbP(\Yhat=\infty)$, where $Z^{(m)}_{t_m}$ is the number of infectious households in generation $t_m$ of $E^{(m)}$. It follows that, with probability tending to $1$ as $m\to\infty$, a major outbreak has at least $\log m$ and at most $(\log m)^\beta$ infectious households in generation $t_m$.

We next consider the probability that a typical individual, $i^*$ say, that is susceptible at time $t_m$ in $E^{(m)}$ ultimately becomes infected. We do this by stopping the construction of $E^{(m)}$ at time $t_m$, leaving the $Z^{(m)}_{t_m}$ infectious (`live') half-edges unconnected, and constructing the susceptibility set, $\calS^{(m)}$ say, of $i^*$ in `generations' as described in Section~\ref{sec:MainBwd}, pairing up the half-edges as we construct the susceptibility set. If at any point in the construction of $\calS^{(m)}$ a half-edge is paired up with one of the $Z^{(m)}_{t_m}$ live half-edges from the epidemic then $i^*$ is ultimately infected, otherwise $i^*$ is not infected by the epidemic. Note that for any individual, $i$ say, in $\calS^{(m)}$ we need to explore all of $i$'s global neighbours (and not just those that join $\calS^{(m)}$), since if \emph{any} half-edge emanating from $i$ is paired with one of the $Z^{(m)}_{t_m}$ live half-edges then $i^*$ is ultimately infected. Thus we need to construct simultaneously $\calA^{(m)}$, the set of global neighbours of $\calS^{(m)}$, also on a generation basis.

Let $(S^{(m)},A^{(m)})=((S^{(m)}_k, A^{(m)}_k),\, k=0,1,\ldots)$ describe the number of households in successive generations of $(\calS^{(m)}, \calA^{(m)})$. In Section~\ref{sec:construction}, we construct realisations of $(S^{(m)},A^{(m)})$ and an approximating two-type branching process $(X^{(m)},X^{(m)}_A)=((X^{(m)}_k,X^{(m)}_{Ak}),\, k=0,1,\ldots)$. The process $X^{(m)}$ is a single-type branching process that is similar to the branching process, $X$ say, described in Section~\ref{sec:MainBwd}, except, as with $Y^{(m)}$, the empirical distribution of $D_1,D_2,\ldots,D_{mn}$ is used instead of the degree distribution $D$. The process $X_A^{(m)}$ corresponds to global neighbours of $\calS^{(m)}$ who are not in $\calS^{(m)}$; individuals in $X_A^{(m)}$ have no offspring. The processes $(S^{(m)},A^{(m)})$ and $(X^{(m)},X_A^{(m)})$ are coupled so that they coincide until $\taubar^{(m)}+1$ households have joined $\calS^{(m)} \cup \calA^{(m)}$, where $\bbP(\taubar^{(m)}>k) \to 1$ as $m\to\infty$ for all $k\in\bbZ_+$. Let $\What^{(m)}$ and $\What_A^{(m)}$ denote the number of households in $\calS^{(m)}$ and $\calA^{(m)}$, respectively, and let $\Xhat^{(m)}$ and $\Xhat$ denote the total progenies of $X^{(m)}$ and $X$, respectively. As with $E^{(m)}$, $\What^{(m)}$ and $\Xhat^{(m)}$ have the same limiting distribution as $m\to\infty$, which, since $X^{(m)}$ converges in distribution to $X$ as $m\to\infty$, is given by the distribution of $\Xhat$. Now, for any $k\in\bbN$, if $\What^{(m)} + \What_A^{(m)} \leq k$ then the probability that $\calS^{(m)}$ intersects with one of the $Z^{(m)}_{t_m}$ live half-edges tends to $0$ as $m\to\infty$ (since a major outbreak has at most $(\log m)^\beta$ infectious households at generation $t_m$ of the forward process), so the limiting (as $m\to\infty$) probability that $i^*$ is ultimately infected by a major outbreak is at most $\bbP(\Xhat=\infty)$. (Note that $(X^{(m)},X_A^{(m)})$ goes extinct if and only if $X^{(m)}$ goes extinct.)

We also construct, for all sufficiently small $\epsilon \in (0,1)$, a branching process ${}_{\epsilon}X^{(m)}$, which is a lower bound for $S^{(m)}$ as long as $\What^{(m)} \leq \epsilon m$; whence $\bbP(\What^{(m)} > \epsilon m) \geq \bbP({}_{\epsilon}\Xhat^{(m)}=\infty)$, where ${}_{\epsilon}\Xhat^{(m)}$ denotes the total progeny of ${}_{\epsilon}X^{(m)}$. As $m\to\infty$, ${}_{\epsilon}\Xhat^{(m)}$ converges in distribution to ${}_{\epsilon}\Xhat$, the total progeny of a branching process ${}_{\epsilon}X$ say. Moreover, for any $\epsilon>0$, if $\What^{(m)} >\epsilon m$ the probability that $\calS^{(m)}$ intersects one of the $Z^{(m)}_{t_m}$ live half-edges tends to $1$ as $m\to\infty$ (since a major outbreak has at least $\log m$ infectious households at generation $t_m$ of the forward process), so the limiting probability that $i^*$ is ultimately infected is at least $\bbP({}_{\epsilon}\Xhat=\infty)$. Furthermore, $\bbP({}_{\epsilon}\Xhat=\infty) \to \bbP(\Xhat=\infty)$ as $\epsilon\downarrow0$, which, combined with the result described at the end of the previous paragraph, shows that the probability that $i^*$ is ultimately infected by a major outbreak tends to $\bbP(\Xhat=\infty)$ as $m\to\infty$ (see Theorem~\ref{thm:majOBInfProb}). It follows that the expected proportion of the population that are infected by a major outbreak also tends to $\bbP(\Xhat=\infty)$ as $m\to\infty$ (see Corollary~\ref{cor:UncondMajOBInfProb}).

Our results are proved by conditioning on the degree sequence $\bD$ and showing that they hold for $\bbP$-almost all $\bD$. The unconditional results then follow using the dominated convergence theorem. As remarked above, the network of global contacts is now constructed as the epidemic/susceptibility set evolves, not a priori as in our model description in Section~\ref{sec:Model}. This implicitly means that, rather than conditioning on the total number of half-edges $\sum_{i=1}^{mn} D_i$ being even, we simply ignore the single left-over half-edge in the event of $\sum_{i=1}^{mn} D_i$ being odd. This small change does not affect the asymptotic results as $m\to\infty$ (cf.~van der Hofstad \etal\ (2007, Section~1.1)).

The remainder of this section is organised as follows. The main constructions are described in Section~\ref{sec:construction}, with the epidemics $E^{(m)}$ and their approximating branching processes being described in Section~\ref{sec:fwdConstruction} and the susceptibility set processes and their approximating branching processes being described in Section~\ref{sec:bwdConstruction}. Some notation concerning the offspring distributions of various branching processes is given in Section~\ref{sec:notAndLimProc}. Section~\ref{sec:PrelimResults} contains some preliminary results, and the main results are given in Sections~\ref{sec:fwdAnalysis} and~\ref{sec:bwdAnalysis}, which analyse the epidemics $E^{(m)}$ and the susceptibility set processes $(S^{(m)},A^{(m)})$, respectively.

\subsection{Construction of approximating branching processes}
\label{sec:construction}
Let $(\Omega_1,\calF_1,\bbP_1)$ be a probability space, on which is defined a sequence $\bD=(D_1,D_2,\ldots)$ of independent random variables, each distributed according to the degree distribution $D$. Also let $(\Omega_2,\calF_2,\bbP_2)$ be a probability space, on which are defined the following mutually independent random quantities:
\begin{enumerate}
\item[(i)] for every $(\bd,j)=((d_1,d_2,\ldots,d_n),j) \in \bbZ_+^n \times \{1,2,\ldots,n\}$, a sequence of random variables $\Phi^{(\bds,j)}_1,\Phi^{(\bds,j)}_2,\ldots$, which are independent copies of the random variable $\Phi^{(\bds,j)}$ defined below.
\item[(ii)] for every $(\bd,j) \in \bbZ_+^n \times \{1,2,\ldots,n\}$, a sequence of random variables $(\Psi^{(\bds,j)}_1,\Psi^{(\bds,j)}_{A1}), $ $ (\Psi^{(\bds,j)}_2,\Psi^{(\bds,j)}_{A2}), \ldots$, which are independent copies of the random variable $(\Psi^{(\bds,j)},\Psi^{(\bds,j)}_{A})$ also defined below.
\end{enumerate}
We also require other random variables defined on $(\Omega_2,\calF_2,\bbP_2)$, but these are described only informally because the detail is unnecessary for our proofs.

The random variable $\Phi^{(\bds,j)}$ describes the number of global neighbours with which infectious contact is made by members of a household of individuals with degrees given by $\bd=(d_1,d_2,\ldots,d_n)$ in which individual $j$ is initially infected and is defined as follows.  Let $G$ be the random directed graph on the vertices $V=\{1,2,\ldots,n\}$ obtained as follows. For each vertex $i$ we take an independent realisation, $I_i$ say, of the infectious period distribution $I$ and then put an arc from $i$ to each other vertex in $V$ independently with probability $1-\e^{-\lambda_L I_i}$. Given $G$, let $C_1,C_2,\ldots,C_n$ be independent random variables with $C_i \cond I_1,I_2,\ldots,I_n \sim \mbox{Bin}(d_i',1-\e^{\lambda_G I_i})$, where $d_i'=d_i$ if $i \ne j$ and $d_j'=d_j-1$. Then $\Phi^{(\bds,j)}=\sum_{i=1}^n \ind{j\leadsto i} C_i$, where $j\leadsto i$ denotes the event that there is a path from vertex $j$ to vertex $i$ in $G$ (with the convention that $i\leadsto i$).

In a similar manner, the two components of the random variable $(\Psi^{(\bds,j)},\Psi_A^{(\bds,j)})$ describe the number of global neighbours of the local susceptibility set of individual $j$ in a household of individuals with degrees given by $\bd=(d_1,d_2,\ldots,d_n)$ that do and do not make global infectious contact with their neighbour in that susceptibility set. To this end, let $G$ be the random graph described above and, conditional on $G$, let $B_1, B_2,\ldots,B_n$ be independent random variables with $B_i \sim \mbox{Bin}(d_i',p_G)$, where $d_1',d_2',\dots,d_n'$ are as above and $p_G=1-\phi(\lambda_G)$. We then have $(\Psi^{(\bds,j)},\Psi_A^{(\bds,j)}) = \sum_{i=1}^n \ind{i\leadsto j} (B_i,d'_i-B_i)$.

We now introduce some further notation. For $k=1,2,\ldots$, let $\bD_k = (D_{k1},D_{k2},\ldots,D_{kn})$, where, for $i=1,2,\ldots,n$, $D_{ki}=D_{(k-1)n+i}$ is the degree of the $i$th individual in the $k$th household. Let $H_k = \sum_{i=1}^n D_{ki}$ denote the total degree of the $k$th household. Lastly, denote by $\mu_{H}^{(m)}= \frac{1}{m} \sum_{i=1}^m H_i$ the (empirical) mean number of edges emanating from each of the first $m$ households.

The epidemic, susceptibility sets and approximating branching processes are defined on the probability space $(\Omega,\calF,\bbP)=(\Omega_1,\calF_1,\bbP_1) \times (\Omega_2,\calF_2,\bbP_2)$. Our construction and most of our calculations will henceforth be conditional on the degree sequence $\bD$. To this end, we denote $\bbP(\cdot \cond \bD)$ by $\bbPD(\cdot)$ and similarly $\bbE[\cdot \cond \bD]$ by $\bbED[\cdot]$. Conditional on this degree sequence and for every $m=1,2,\ldots$, we now describe the construction of a branching process, $Y^{(m)}$, which approximates the early stages of the spread of the epidemic amongst households $1,2,\ldots,m$; then another (two-type) branching process, $(X^{(m)},X_A^{(m)})$, which approximates the `early growth' of the susceptibility set (and its global neighbours) of a typical initially susceptible individual in that population.

\subsubsection{The forward processes.}
\label{sec:fwdConstruction}
We first describe the branching process $Y^{(m)}$. Set $Y^{(m)}_0=1$ and choose an individual uniformly at random from $1,2,\ldots,mn$. Suppose it is individual $\iota\in\{1,2,\ldots,n\}$ of household $\Delta_0\in\{1,2,\ldots,m\}$. Then $Y^{(m)}_1=\Phi^{(\bDs_{\Delta_0}+e_\iota,\iota)}_1$, where $e_i$ is the unit $n$-vector with a $1$ in the $i$th position. For subsequent generations $k\geq2$, we continue the construction as follows. For each $j=1,2,\ldots,Y^{(m)}_{k-1}$, sample a half-edge uniformly at random from the $m\mu_{H}^{(m)}$ half-edges in the population and, supposing it emanates from individual $\iota$ of household $\Delta$, set $Y^{(m)}_{kj}=\Phi^{(\bDs_{\Delta},\iota)}_{\nu(\Delta,\iota)+1}$, where $\nu(\Delta,\iota)$ is the number of times we have sampled previously from the sequence $\Phi^{(\bDs_{\Delta},\iota)}_1,\Phi^{(\bDs_{\Delta},\iota)}_2,\ldots$. Lastly, set
\begin{equation*}
Y^{(m)}_k = \sum_{j=1}^{Y^{(m)}_{k-1}} Y^{(m)}_{kj}.
\end{equation*}

The branching process $Y^{(m)}$ and the epidemic process $E^{(m)}$ can be coupled by using the same $\bD$, $\Phi$'s and uniformly random samples. However, the coupling breaks down as soon as a half-edge is sampled that emanates from a household that either has been used previously in the epidemic or is a neighbour of such a previously used household. If a previously used half-edge is sampled then in $E^{(m)}$ another half-edge needs to be sampled. If an unused half-edge that emanates from a previously used household is sampled then in $E^{(m)}$ the spread of the epidemic within that household is different from in $Y^{(m)}$ since there are fewer susceptibles. Finally, if a half-edge emanating from a household neighbouring a household previously used in $E^{(m)}$ is sampled then the spread of the epidemic from that household is in general different from that in $Y^{(m)}$, since the (effective) degree distribution of individuals in that household may be different from that assumed in $Y^{(m)}$. (When constructing $E^{(m)}$ one needs also to pair up non-infectious half-edges from infectious individuals.) In all of these cases the construction of $E^{(m)}$ can be continued appropriately but the detail is not important for our purposes. However, we do need a bound on the size of, and number of half-edges that emanate from, the `bad set' of households that must be avoided in order that $Y^{(m)}$ and $E^{(m)}$ remain coupled. To that end we describe another branching process $T^{(m)}=(T^{(m)}_k,\, k=1,2,\ldots)$, which provides such a bound.

Let $T^{(m)}_0$ be the total degree of the initial household in $Y^{(m)}$, so $T^{(m)}_0 = H_{\Delta_0}$, where $\Delta_0$ is as above. For $k=1,2,\ldots$, $T^{(m)}_k$ is determined as follows. For each $j=1,2,\ldots,T^{(m)}_{k-1}$, a half-edge is sampled uniformly at random from the $m\mu_{H}^{(m)}$ half-edges in the population, say this half-edge emanates from household $\Delta_j$, and then put $T^{(m)}_{kj}=H_{\Delta_j}-1$. Finally, set $T^{(m)}_k = \sum_{j=1}^{T^{(m)}_{k-1}} T^{(m)}_{kj}$. The processes $Y^{(m)}$, $E^{(m)}$ and $T^{(m)}$ can be coupled in an obvious fashion so that their sampled half-edges correspond. Let $\That^{(m)}_k = \sum_{l=0}^k T^{(m)}_l$ ($k=0,1,\ldots$) be the total progeny of $T^{(m)}$ up to generation $k$. Then $2\That^{(m)}_{k+1}$ provides an upper bound for the number of half-edges that emanate from (and hence also for the size) of the bad set of households in generation $k$ of $E^{(m)}$. The index $k+1$ arises because the bad set consists of not just all households infected up to generation $k$ of $E^{(m)}$ but also their neighbouring households. The factor $2$ arises because $T^{(m)}$ does not count the receiving half-edge when the half-edges are paired up.

The above construction of $Y^{(m)}$ (and implicitly $E^{(m)}$) is continued for a fixed number of generations, $t_m$, and $T^{(m)}$ is continued for $t_m+1$ generations. (Of course, some or all of these processes may die out beforehand.)

\subsubsection{The backward processes.}
\label{sec:bwdConstruction}
The two-type branching process $(X^{(m)},X_A^{(m)})$ is defined analogously to $Y^{(m)}$ except the random variables $(\Psi^{(\bds,j)}_i,\Psi^{(\bds,j)}_{Ai})$ are used instead of $\Phi^{(\bds,j)}_i$ (recall that there are no offspring in $X_A^{(m)}$). The process $X^{(m)}$ approximates the growth, described by generations as in Section~\ref{sec:MainBwd}, of the susceptibility set of an individual chosen uniformly at random from all susceptible individuals at time $t_m$ in the epidemic process $E^{(m)}$ and $X_A^{(m)}$ approximates the number of global neighbours of this susceptibility set, also on a generation basis. The processes $(X^{(m)},X_A^{(m)})$ and $(S^{(m)},A^{(m)})$ can be coupled in a similar fashion to that used for $Y^{(m)}$ and $E^{(m)}$, though note that now the coupling breaks down if a sampled half-edge emanates from either (i) a household previously used in the susceptibility set, (ii) a household neighbouring such a household, or (iii) a household or neighbour of a household used in the forward process up to time $t_m$. Also note that this coupling may break down at generation $0$ (if the initial individual is in a household that is either infected in $E^{(m)}$ or a neighbour of a household infected in $E^{(m)}$). As with the epidemic process $E^{(m)}$, the construction of $(S^{(m)},A^{(m)})$ can be continued appropriately after the coupling breaks down but we do not require such detail.

\subsubsection{Further notation and limiting processes.}
\label{sec:notAndLimProc}
For $m=1,2,\ldots$, let $\bc^{(m)}=(c^{(m)}_0,c^{(m)}_1,\ldots)$ and $\bctilde^{(m)}=(\ctilde^{(m)}_0,\ctilde^{(m)}_1,\ldots)$ denote the offspring distributions of the initial individual and all subsequent individuals, respectively, in $Y^{(m)}$. For $\bd=(d_1,d_2,\ldots,d_n)\in\bbZ_+^n$, let
\begin{equation*}
p_{\bds}^{(m)} = \frac{1}{m} \sum_{i=1}^m \ind{\bDs_i=\bds}
\qquad \mbox{and} \qquad
\ptilde_{\bds}^{(m)} = \frac{\abs{\bd}}{m\mu_{H}^{(m)}} \sum_{i=1}^m \ind{\bDs_i=\bds},
\end{equation*}
where $\abs{\bd}=\sum_{j=1}^n d_j$. Then the `household type' (i.e.~the degrees of individuals within the household) of the initial individual in $Y^{(m)}$ is distributed according to $p_{\bds}^{(m)}$ ($\bd\in\bbZ_+^n$) and the household type of any subsequent individual is distributed according to $\ptilde_{\bds}^{(m)}$ ($\bd\in\bbZ_+^n$). It follows that, for $k=0,1,\ldots$,
\begin{equation}
c^{(m)}_k = \sum_{\bds\in\bbZ_+^n} p_{\bds}^{(m)} \bbP(\Phi_{\bds}=k)
\qquad \mbox{and} \qquad
\ctilde^{(m)}_k = \sum_{\bds\in\bbZ_+^n} \ptilde_{\bds}^{(m)} \bbP(\Phitilde_{\bds}=k),
\label{eq:cmkAndctildemkDef}
\end{equation}
where $\Phi_{\bds}$ and $\Phitilde_{\bds}$ are random variables with distributions given by
\begin{equation}
\bbP(\Phi_{\bds}=k) = \frac{1}{n} \sum_{i=1}^n \bbP(\Phi^{(\bds+e_i,i)}=k) \qquad (k=0,1,\ldots,\abs{\bd}),
\label{eq:PhidDef}
\end{equation}
and
\begin{equation}
\bbP(\Phitilde_{\bds}=k) = \sum_{i=1}^n \frac{d_i}{\abs{\bd}} \bbP(\Phi^{(\bds,i)}=k) \qquad (k=0,1,\ldots,\abs{\bd}-1).
\label{eq:PhitildedDef}
\end{equation}

For $m=1,2,\ldots$, the offspring distributions of the initial and subsequent individuals in $X^{(m)}$, $\bb^{(m)}$ and $\bbtilde^{(m)}$, are defined analagously to $\bc^{(m)}$ and $\bctilde^{(m)}$, using~\eqref{eq:cmkAndctildemkDef}--\eqref{eq:PhitildedDef} with $\Phi$ replaced by $\Psi$ and $\Phitilde$ by $\Psitilde$ throughout. Replacing $\Phi$ by $(\Psi,\Psi_A)$ and $\Phitilde$ by $(\Psitilde,\Psitilde_A)$ throughout gives the offspring distributions associated with the two-type process $(X^{(m)},X_A^{(m)})$.

Further, for $m=1,2,\ldots$, let $\br^{(m)}=(r^{(m)}_0,r^{(m)}_1,\ldots)$ denote the distribution of the number of initial ancestors and $\brtilde^{(m)}=(\rtilde^{(m)}_0,\rtilde^{(m)}_1,\ldots)$ denote the offspring distribution of both the ancestors and any subsequent individuals in $T^{(m)}$. Then, for $k=0,1,\ldots$,
\begin{equation*}
r^{(m)}_k = \sum_{\{\bds\in\bbZ_+^n \,:\, \abs{\bds}=k\}} p_{\bds}^{(m)}
\qquad \mbox{and} \qquad
\rtilde^{(m)}_k = \sum_{\{\bds\in\bbZ_+^n \,:\, \abs{\bds}=k+1\}} \ptilde_{\bds}^{(m)}.
\end{equation*}
For $m=1,2,\ldots$, let $\mu_c^{(m)} = \sum_{k=1}^\infty k c^{(m)}_k$ be the mean of the empirical distribution $\bc^{(m)}$, and define $\mutilde_c^{(m)}, \mu_b^{(m)},\mutilde_b^{(m)},\mu_r^{(m)}$ and $\mutilde_r^{(m)}$ analogously.

In Section~\ref{sec:PrelimResults} we prove that the offspring distributions of $Y^{(m)}$, $X^{(m)}$ and $T^{(m)}$ and the distribution of the number of ancestors in $T^{(m)}$ converge almost surely as $m\to\infty$ to those of branching processes we denote by $Y$, $X$ and $T$ respectively. To that end, for $\bd\in\bbZ_+^n$, let $p_{\bds} = \prod_{i=1}^n p_{d_i}$ and $\ptilde_{\bds} = p_{\bds} \abs{\bd}/n\mu_D$. (Recall that $p_k=\bbP(D=k)$ ($k=0,1,\ldots$) and $\mu_D=\sum_{k=1}^\infty kp_k$.) Also, for $k=0,1,\ldots$, let
\begin{equation*}
p_H(k) = \sum_{\{\bds\in\bbZ_+^n \,:\, \abs{\bds}=k\}} p_{\bds} = \bbP(D_1+D_2+\cdots+D_n=k) =\bbP(H_1=k)
\end{equation*}
and, for $k=1,2,\ldots$, let $\ptilde_H(k) = kp_H(k)/n\mu_D$. Now, for $k=0,1,\ldots$, let $c_k$ be defined analogously to $c^{(m)}_k$ but with $p_{\bds}^{(m)}$ replaced by $p_{\bds}$, and define $\ctilde_k, b_k$ and $\btilde_k$ similarly. Also, for $k=0,1,\ldots$, let $r_k = p_H(k)$ and $\rtilde_k = \ptilde_H(k+1)$. Let $\bc=(c_0,c_1,\ldots)$ and define $\bctilde$, $\bb$, $\bbtilde$, $\br$ and $\brtilde$ similarly. Let $\mu_c = \sum_{k=1}^\infty k c_k$ and define $\mutilde_c$, $\mu_b$, $\mutilde_b$, $\mu_r$ and $\mutilde_r$ in the obvious fashion.

Let $Y=(Y_0,Y_1,\ldots)$, $X=(X_0,X_1,\ldots)$ and $T=(T_0,T_1,\ldots)$ be the branching processes ${\rm BP}(1,\bc,\bctilde)$, ${\rm BP}(1,\bb,\bbtilde)$ and, in an obvious notation, ${\rm BP}(\br,\brtilde,\brtilde)$, respectively. Note that the branching processes $Y$ and $X$ are those described in Sections~\ref{sec:MainFwd} and~\ref{sec:MainBwd}, respectively. Note especially that, in the notation of Section~\ref{sec:MainFwd}, this implies that $\mutilde_c=R_*$. We also require a two-type branching process $(X,A)$, defined analagously to $(X^{(m)},X_A^{(m)})$ but again using $p_{\bds}$ and $\ptilde_{\bds}$ in defining the offspring distribution instead of the empirical versions $p_{\bds}^{(m)}$ and $\ptilde_{\bds}^{(m)}$.

\subsection{Preliminary results}
\label{sec:PrelimResults}
In this section we collect some results required in the analysis of the forward and backward processes. Recall that we have made the assumption that $\sigma^2_D=\Var D$ is finite (although some results only require $\mu_D<\infty$).
\begin{lem}
\label{lem:SLLN}
There exists $A_1\in\calF_1$, with $\bbP_1(A_1)=1$, such that, for all $\omega_1\in A_1$,
\begin{enumerate}
\item[(i)] $\ds \lim_{m\to\infty} \mu_{H}^{(m)}(\omega_1) = n\mu_D$;
\item[(ii)] $\ds \lim_{m\to\infty} p_{\bds}^{(m)}(\omega_1) = p_{\bds}$ and $\ds \lim_{m\to\infty} \ptilde_{\bds}^{(m)}(\omega_1) = \ptilde_{\bds}$ for each $\bd\in\bbZ_+^n$;
\item[(iii)] $\ds \lim_{m\to\infty} \bc^{(m)}(\omega_1) = \bc$, $\ds \lim_{m\to\infty} \bctilde^{(m)}(\omega_1) = \bctilde$, $\ds \lim_{m\to\infty} \bb^{(m)}(\omega_1) = \bb$, \\$\ds \lim_{m\to\infty} \bbtilde^{(m)}(\omega_1) = \bbtilde$, $\ds \lim_{m\to\infty} \br^{(m)}(\omega_1) = \br$ and $\ds \lim_{m\to\infty} \brtilde^{(m)}(\omega_1) = \brtilde$;
\item[(iv)] $\ds \lim_{m\to\infty} \mu_c^{(m)}(\omega_1) = \mu_c$, $\ds \lim_{m\to\infty} \mutilde_c^{(m)}(\omega_1) = \mutilde_c$, $\ds \lim_{m\to\infty} \mu_b^{(m)}(\omega_1) = \mu_b$, \\$\ds \lim_{m\to\infty} \mutilde_b^{(m)}(\omega_1) = \mutilde_b$, $\ds \lim_{m\to\infty} \mu_r^{(m)}(\omega_1) = \mu_r$ and $\ds \lim_{m\to\infty} \mutilde_r^{(m)}(\omega_1) = \mutilde_r$.
\end{enumerate}
Here and henceforth, convergence of a sequence of sequences is interpreted elementwise, so, for example, $\lim_{m\to\infty} \bc^{(m)} = \bc$ means that $\lim_{m\to\infty} c^{(m)}_k = c_k$ for each $k=0,1,\ldots$.
\end{lem}
\begin{proof}
By the strong law of large numbers, there exists $A_2\in\calF_1$ with $\bbP_1(A_2)=1$ such that $\ds \lim_{m\to\infty} \mu_{H}^{(m)}(\omega_1) = n\mu_D$ ($\omega_1\in A_2$) and, for each $\bd\in\bbZ_+^n$, there exists $A_{\bds}\in\calF_1$ with $\bbP_1(A_{\bds})=1$ such that $\ds \lim_{m\to\infty} p_{\bds}^{(m)}(\omega_1) = p_{\bds}$ ($\omega_1\in A_{\bds}$). Let $A_3=A_2\cap \bigcap_{\bds\in\bbZ_+^n} A_{\bds}$. Then $\bbP_1(A_3)=1$ and it is easily verified that (i) and (ii) hold for all $\omega_1\in A_3$, whence (iii) also holds for all $\omega_1\in A_3$ by \Scheffe's theorem (see, for example, Billingsley (1968, p.~224)). Next, consider
\begin{eqnarray*}
\mutilde_c^{(m)} = \sum_{k=1}^\infty k \ctilde^{(m)}_k & = & \sum_{k=1}^\infty k \sum_{\bds\in\bbZ_+^n} \ptilde_{\bds}^{(m)} \bbP(\Phitilde_{\bds}=k) \\
 & = & \sum_{k=1}^\infty k \sum_{\bds\in\bbZ_+^n} \frac{\abs{\bd}}{m\mu_{H}^{(m)}} \sum_{i=1}^m \ind{\bDs_i=\bds} \bbP(\Phitilde_{\bds}=k) \\
 & = & \frac{1}{\mu_{H}^{(m)}} \cdot \frac{1}{m} \sum_{i=1}^m \abs{\bD_i} \sum_{k=1}^\infty k \bbP(\Phitilde_{\bDs_i}=k).
\end{eqnarray*}
Now, $\bbP(\Phitilde_{\bDs_i}=k)=0$ for $k\geq \abs{\bD_i}-1$, so
\begin{equation*}
 \bbE\b{ \abs{\bD_1} \sum_{k=1}^\infty k \bbP(\Phitilde_{\bDs_1}=k) } \leq \bbE\b{ \abs{\bD_1}(\abs{\bD_1}-1) } <\infty,
\end{equation*}
as $\sigma^2_D<\infty$. Thus, by the strong law of large numbers, there exists $A_4\in\calF_1$ with $\bbP_1(A_4)=1$ such that, for all $\omega_1\in A_4$,
\begin{eqnarray*}
\lim_{m\to\infty} \frac{1}{m} \sum_{i=1}^m \abs{\bD_i(\omega_1)} \sum_{k=1}^\infty k \bbP(\Phitilde_{\bDs_i(\omega_1)}=k) & = & \bbE \b{ \abs{\bD_1} \sum_{k=1}^\infty k \bbP(\Phitilde_{\bDs_1}=k) } \\
 & = & \sum_{\bds\in\bbZ_+^n} p_{\bds} \abs{\bd} \sum_{k=1}^\infty k \bbP(\Phitilde_{\bds}=k).
\end{eqnarray*}
It follows that $ \lim_{m\to\infty} \mutilde_c^{(m)}(\omega_1) = \mutilde_c$ for all $\omega_1\in A_2\cap A_4$. Similar arguments hold for the other means in (iv) and the lemma is thus proved.
\end{proof}
\begin{remnn}
Throughout the remainder of the paper, $A_1$ refers to a set that satisfies the statement of Lemma~\ref{lem:SLLN}.
\end{remnn}

The following result concerns the convergence of certain quantities associated with a sequence of branching processes when their offspring distributions converge in distribution.
\begin{lem}
\label{lem:BPTPcgce}
Suppose that $\ba$, $\batilde$, $\ba^{(m)}$ and $\batilde^{(m)}$ ($m=1,2,\ldots$) are probability distributions satisfying $\ba^{(m)} \to \ba$ and $\batilde^{(m)} \to \batilde$ as $m\to\infty$. Let $Y^{(m)}\sim {\rm BP}(1,\ba^{(m)},\batilde^{(m)})$ $(m=1,2,\ldots)$ and $Y\sim {\rm BP}(1,\ba,\batilde)$. Then, denoting by $\Yhat^{(m)}$ (respectively $\Yhat$) the total progeny of $Y^{(m)}$ (respectively $Y$),
\begin{enumerate}
\item[(i)] $\ds \lim_{m\to\infty} \bbP(\Yhat^{(m)}=k) = \bbP(\Yhat=k)$ ($k=1,2,\ldots$);
\item[(ii)] $\ds \lim_{m\to\infty} \bbP(\Yhat^{(m)}=\infty) = \bbP(\Yhat=\infty)$, provided $\atilde_1 \neq 1$.
\end{enumerate}
\end{lem}
\begin{proof}
Part (i) follows immediately by considering the sum of the probabilities of the finite number of sample paths of $Y^{(m)}$ with $\Yhat^{(m)}=k$. Part (ii) is a simple extension of Lemma~4.1 of Britton \etal\ (2007).
\end{proof}

\begin{remsnn}\begin{enumerate}
\item The condition in part (ii) of the lemma is in practice only a technical condition which will always hold true. As pointed out by Britton \etal\ (2007), although the case $\atilde_1=1$ really can be an exception (for example if $\atilde^{(m)}_0 = 1-\atilde^{(m)}_1 = 1/m$), such a scenario is, from an applied viewpoint, decidedly pathological.
\item We sometimes use a slight variant of Lemma~\ref{lem:BPTPcgce}, where the branching processes are indexed by $\epsilon\in(0,1)$ and their offspring distributions converge as $\epsilon\downarrow0$. Of course, the analogous results hold, and the proof is exactly the same.
\end{enumerate}
\end{remsnn}

Lastly, we have a result concerning the probability of picking a `bad' half-edge in our constructions of the forward and backward processes.
\begin{lem}
\label{lem:gh}
Suppose that, for each $m=1,2,\ldots$, we draw elements uniformly at random, with replacement, from the set $\calJ^{(m)}=\{1,2,\ldots,m\mu_{H}^{(m)}\}$. Suppose also that, for each $m$, there is an increasing sequence of (random) sets $\calJ^{(m)}_1 \subset \calJ^{(m)}_2 \subset \cdots \subset \calJ^{(m)}$ and at the $i$th pick we wish to avoid picking a member of $\calJ^{(m)}_i$. Denote the $i$th pick by $\chi^{(m)}_i$ and let $\tau^{(m)}=\min\{i \,:\, \chi^{(m)}_i \in \calJ^{(m)}_i\}-1$ be the number of picks we make before making a pick from a set we wish to avoid. Suppose further that there exist strictly positive integers $g(m)$ and $h(m)$ ($m=1,2,\ldots$) satisfying $\lim_{m\to\infty} g(m)h(m)m^{-1} = 0$ and
\begin{equation}
\lim_{m\to\infty} \bbPDo(J^{(m)}_{g(m)} \leq h(m)) = 1
\label{eq:sizeALeqH}
\end{equation}
for all $\omega_1\in A_1$, where $J^{(m)}_i=|\calJ^{(m)}_i|$. Then, for all $\omega_1\in A_1$,
\begin{equation}
\lim_{m\to\infty} \bbPDo(\tau^{(m)}>g(m)) = 1.
\label{eq:tau>g}
\end{equation}
\end{lem}
\begin{proof}
In view of~\eqref{eq:sizeALeqH}, for $\omega_1\in A_1$,
\begin{align*}
\liminf_{m\to\infty} & \bbPDo(\tau^{(m)}>g(m)) \\
 & = \liminf_{m\to\infty} \bbPDo \pp{ \tau^{(m)}>g(m) \, \big| \, J^{(m)}_{g(m)} \leq h(m) } \bbPDo\pp{J^{(m)}_{g(m)} \leq h(m)} \\
 & \geq \liminf_{m\to\infty} \pp{1-\frac{h(m)}{m\mu_{H}^{(m)}(\omega_1)}}^{g(m)} \bbPDo\pp{J^{(m)}_{g(m)} \leq h(m)} \\
 & \geq \liminf_{m\to\infty} \pp{1-\frac{g(m)h(m)}{m\mu_{H}^{(m)}(\omega_1)}} \bbPDo\pp{J^{(m)}_{g(m)} \leq h(m)} \\
 & = 1,
\end{align*}
using~\eqref{eq:sizeALeqH}, Lemma~\ref{lem:SLLN}(i) and the fact that $g(m)h(m)m^{-1} \to 0$ as $m\to\infty$. The assertion~\eqref{eq:tau>g} then follows.
\end{proof}

\subsection{Analysis of forward process}
\label{sec:fwdAnalysis}
\subsubsection{Threshold theorem for the epidemic $E^{(m)}$.}
\label{sec:threshold}
In order to prove a threshold theorem for the epidemic we first establish a bound for the size of the bad set of half-edges after $k$ generations of the epidemic $E^{(m)}$. Recall (from the discussion at the end of Section~\ref{sec:fwdConstruction}) that the number of half-edges in this set is bounded by $2\That^{(m)}_{k+1}$.
\begin{lem}
\label{lem:ThatLogm}
For all $\omega_1\in A_1$,
\begin{equation*}
\lim_{m\to\infty} \bbPDo(\That^{(m)}_{k} > \log m) =0 \qquad (k=1,2,\ldots).
\end{equation*}
\end{lem}
\begin{proof}
Fix $\omega_1\in A_1$. Then note that $\bbEDo[\That^{(m)}_{k}] = \mu_r^{(m)}(\omega_1)\{1+\mutilde_r^{(m)}(\omega_1) + (\mutilde_r^{(m)}(\omega_1))^2 + \cdots + (\mutilde_r^{(m)}(\omega_1))^{k}\}$ and also that $\mu_r^{(m)}(\omega_1) \leq \mutilde_r^{(m)}(\omega_1)+1$. Thus $\bbEDo[\That^{(m)}_{k}] \leq (k+1) (\mutilde_r^{(m)}(\omega_1) + 1)^{k+1}$ and, by Markov's inequality,
\begin{equation*}
\bbPDo(\That^{(m)}_{k}(\omega_1) > \log m) \leq \frac{k+1}{\log m} (\mutilde_r^{(m)}(\omega_1) + 1)^{k+1}.
\end{equation*}
The lemma now follows, since $\mutilde_r^{(m)}(\omega_1)\to \mutilde_r$ as $m\to\infty$ for all $\omega_1\in A_1$, by Lemma~\ref{lem:SLLN}(iv).
\end{proof}

For $m=1,2,\ldots$, let $\Zhat^{(m)}$ denote the total number of households infected in the epidemic $E^{(m)}$, including the initial household, and let $\Yhat^{(m)}$ and $\Yhat$ be the total progeny, including the initial individual, of the branching processes $Y^{(m)}$ and $Y$, respectively. We now show that the total number of households infected in $E^{(m)}$ converges in distribution to the total progeny of $Y$.
\begin{thm}
\label{thm:threshold}
For $k=1,2,\ldots$,
\begin{enumerate}
\item[(i)] for all $\omega_1\in A_1$, $\ds \lim_{m\to\infty} \bbPDo(\Zhat^{(m)}=k) = \bbP(\Yhat=k)$;
\item[(ii)] $\ds \lim_{m\to\infty} \bbP(\Zhat^{(m)}=k) = \bbP(\Yhat=k)$.
\end{enumerate}
\end{thm}
\begin{proof}
Fix $\omega_1\in A_1$ and let $\tau^{(m)}$ be the number of households infected by $E^{(m)}$ before a bad half-edge is chosen. Fix $k\in\bbN$. Then
\begin{equation}
\bbPDo(\Zhat^{(m)}=k) = \bbPDo(\Zhat^{(m)}=k \comma \tau^{(m)}\leq k) + \bbPDo(\Zhat^{(m)}=k \comma \tau^{(m)}> k).
\label{eq:splitPDoZhatm}
\end{equation}
Let $\calJ^{(m)}_l$ ($l=1,2,\ldots$) be the set of half-edges we wish to avoid when choosing the $l$th household to spread the epidemic to. Then $J^{(m)}_k=|\calJ^{(m)}_k| \leq 2\That^{(m)}_{k}$, so
\begin{equation*}
\bbPDo(J^{(m)}_{k} \leq 2\log m) \geq \bbPDo(\That^{(m)}_{k} \leq \log m) \to 1
\end{equation*}
as $m\to\infty$, by Lemma~\ref{lem:ThatLogm}. Thus, using Lemma~\ref{lem:gh} with $g(m)=k$ and $h(m)=2\log m$, $\lim_{m\to\infty} \bbPDo(\tau^{(m)}> k) = 1$. Therefore, $\lim_{m\to\infty} \bbPDo(\Zhat^{(m)}=k \comma \tau^{(m)}\leq k) = 0$ and, recalling~\eqref{eq:splitPDoZhatm},
\begin{eqnarray*}
\lim_{m\to\infty} \bbPDo(\Zhat^{(m)}=k) & = & \lim_{m\to\infty} \bbPDo(\Zhat^{(m)}=k \comma \tau^{(m)}> k) \\
 & = & \lim_{m\to\infty} \bbPDo(\Yhat^{(m)}=k \comma \tau^{(m)}> k) \\
 & = & \lim_{m\to\infty} \bbPDo(\Yhat^{(m)}=k) \\
 & = & \bbP(\Yhat=k),
\end{eqnarray*}
using Lemmas~\ref{lem:SLLN}(iii) and~\ref{lem:BPTPcgce}(i), proving assertion (i). Further,
\begin{equation*}
\lim_{m\to\infty} \bbP(\Zhat^{(m)}=k) = \lim_{m\to\infty} \bbE\b{ \bbPD(\Zhat^{(m)}=k) } = \bbP(\Yhat=k),
\end{equation*}
using the dominated convergence theorem, proving assertion (ii).
\end{proof}

\subsubsection{Early behaviour of major outbreaks.}
\label{sec:earlyMajOB}
Theorem~\ref{thm:threshold} shows that the total number of households infected in $E^{(m)}$ converges in distribution as $m\to\infty$ to the total progeny of $Y$, so if $\mutilde_c\leq1$ only minor outbreaks can occur in the limit as $m\to\infty$ (recall that $R_*=\mutilde_c$). We now assume that $\mutilde_c>1$ and study the early behaviour of $E^{(m)}$ when a major outbreak occurs. For $m=1,2,\ldots$, let
\begin{equation*}
t_m=\floor{2\log \log m/\log \mutilde_c}.
\end{equation*}
We obtain a bound on the size of the `bad set' of half-edges at time $t_m$ and show that, with probability tending to $1$ as $m\to\infty$, in a major outbreak there are at least $\log m$ infected households after $t_m$ generations of the epidemic process.
\begin{lem}
\label{lem:ThatLogm_beta}
There exists $\beta\in(1,\infty)$ such that, for all $\omega_1\in A_1$,
\begin{equation*}
\lim_{m\to\infty} \bbPDo(\That^{(m)}_{t_m+1} \geq (\log m)^\beta) = 0.
\end{equation*}
\end{lem}
\begin{proof}
Fix $\omega_1\in A_1$ and note that, for all sufficiently large $m$,
\begin{eqnarray*}
\bbEDo[\That^{(m)}_{t_m+1}] & = & \mu_r^{(m)}(\omega_1)(1+\mutilde_r^{(m)}(\omega_1) + \cdots + (\mutilde_r^{(m)}(\omega_1))^{t_m+1}) \\
 & \leq & \mu_r^{(m)}(\omega_1) \frac{(\mutilde_r^{(m)}(\omega_1))^{t_m+2}}{\mutilde_r^{(m)}(\omega_1)-1}.
\end{eqnarray*}
Thus, by Markov's inequality, for such $m$,
\begin{equation}
\label{eq:logmbeta}
\bbPDo(\That^{(m)}_{t_m+1} \geq (\log m)^\beta) \leq \frac{(\mutilde_r^{(m)}(\omega_1))^{t_m+2}}{(\log m)^\beta} \, \frac{\mu_r^{(m)}(\omega_1)}{\mutilde_r^{(m)}(\omega_1)-1}.
\end{equation}
It is readily shown, by considering its logarithm and using Lemma~\ref{lem:SLLN}(iv), that, for all sufficiently large $\beta$, the right hand side of \eqref{eq:logmbeta} tends to $0$ as $m \to \infty$, and the lemma follows.
\end{proof}

\begin{lem}
\label{lem:limProbYt}
For all $\omega_1\in A_1$,
\begin{equation*}
\lim_{m\to\infty} \bbPDo(Y^{(m)}_{t_m} > \log m) = \bbP(\Yhat=\infty).
\end{equation*}
\end{lem}
\begin{proof}
Note that either (i) $\bc$ and $\bctilde$ both have infinite support, or (ii) $\bc$ and $\bctilde$ are supported on $\{0,1,\ldots,n\dmax\}$ and $\{0,1,\ldots,n\dmax-1\}$ (or subsets thereof) respectively, where $\dmax=\max\{k \,:\, p_k>0\}$.

Consider (i) first. For sufficiently small $\epsilon>0$, let $k_0=\min\{k \,:\, \sum_{i=k+1}^\infty c_i<\epsilon\}$, $\epsilon'=\sum_{i=k_0+1}^\infty c_i$, $\ktilde_0=\min\{k \,:\, \sum_{i=k+1}^\infty \ctilde_i<\epsilon\}$ and $\epsilontilde'=\sum_{i=\ktilde_0+1}^\infty \ctilde_i$, so $\epsilon' < \epsilon$ and $\epsilontilde' < \epsilon$. Note that $k_0$ and $\ktilde_0$ are well-defined whenever $\epsilon<1-(c_0\vee\ctilde_0)$ (where $a \vee b = \max(a,b)$), and also that both $k_0$ and $\ktilde_0$ tend to $\infty$ as $\epsilon\downarrow 0$. Now let $Y^\epsilon=(Y^\epsilon_k,\, k=0,1,\ldots)\sim {\rm BP}(1,\bc^{\epsilon},\bctilde^{\epsilon})$, where $\bc^\epsilon$ has elements $c_i^\epsilon=c_i+\frac{\epsilon'}{k_0+1}$ for $i=0,1,\ldots,k_0$ and $c_i^\epsilon=0$ for $i>k_0$, and $\bctilde^\epsilon=(\ctilde_i^\epsilon,\,i=0,1,\ldots)$ is defined similarly but with $c_i, \epsilon'$ and $k_0$ replaced by $\ctilde_i,\epsilontilde'$ and $\ktilde_0$, respectively. Also let $\mu_\epsilon=\sum_{k=1}^\infty kc^\epsilon_k$ and $\mutilde_\epsilon=\sum_{k=1}^\infty k\ctilde^\epsilon_k$.

Now, note that $\sum_{i=0}^k c_i < \sum_{i=0}^k c_i^\epsilon$ ($k=0,1,\ldots$), so $\mu_\epsilon<\mu_c$. We also have $\mu_\epsilon = \sum_{i=1}^{k_0} i c_i^\epsilon \geq \sum_{i=1}^{k_0} i c_i \to \mu_c$ as $\epsilon\downarrow0$, so $\mu_\epsilon\to\mu_c$ as $\epsilon\downarrow0$. Similarly, $\mutilde_\epsilon\to\mutilde_c$ as $\epsilon\downarrow0$. Now fix $\omega_1\in A_1$. Then, by Lemma~\ref{lem:SLLN}(iii), $\bc^{(m)}(\omega_1)\to \bc$ and $\bctilde^{(m)}(\omega_1)\to \bctilde$ as $m\to\infty$, so there exists $M(\epsilon,\omega_1)$ such that, for all $m\geq M(\epsilon,\omega_1)$, $c^{(m)}_i(\omega_1) < c_i^\epsilon$, for $i=1,2,\ldots,k_0$, and $\ctilde^{(m)}_i(\omega_1) < \ctilde_i^\epsilon$, for $i=1,2,\ldots,\ktilde_0$. Thus, for $m\geq M(\epsilon,\omega_1)$ and $k=0,1,\ldots$, $\sum_{i=0}^k c^{(m)}_i(\omega_1) < \sum_{i=0}^k c_i^\epsilon$ and $\sum_{i=0}^k \ctilde^{(m)}_i(\omega_1) < \sum_{i=0}^k \ctilde_i^\epsilon$ (note, for example, that $\sum_{i=0}^k c_i^\epsilon = 1$ if $k\geq k_0$), whence $Y^{(m)}(\omega_1) \geqst Y^\epsilon$, where $\geqst$ denotes stochastic ordering. Therefore, for $\omega_1\in A_1$ and $m\geq M(\epsilon,\omega_1)$,
\begin{equation}
\bbPDo(Y^{(m)}_{t_m} > \log m) \geq \bbP(Y^\epsilon_{t_m} > \log m) = \bbP\pp{\frac{Y^\epsilon_{t_m}}{\mu_\epsilon \mutilde_\epsilon^{t_m-1}} > \frac{\log m}{\mu_\epsilon \mutilde_\epsilon^{t_m-1}}}.
\label{eq:YlogUpper}
\end{equation}

Now, note that $\sum_{i=1}^\infty \ctilde^\epsilon_i i\log i <\infty$ (the summand is $0$ for $i>\ktilde_0$), so by the well known result concerning the exponential growth of branching processes (see, for example, Haccou \etal\ (2005, Theorem~6.1)), there exists a random variable $\calW^\epsilon$, which takes the value $0$ if and only if $Y^\epsilon$ goes extinct (i.e.~if and only if $\Yhat^\epsilon=\sum_{i=0}^\infty Y^\epsilon_i <\infty$), such that
\begin{equation*}
\frac{Y^\epsilon_{t_m}}{\mu_\epsilon \mutilde_\epsilon^{t_m-1}} \toas \calW^\epsilon \qquad \mbox{as $m\to\infty$}.
\end{equation*}
Next, since $t_m = \floor{2\log\log m / \log \mutilde_c}$, observe that, for suitable $\theta_m\in[0,1)$,
\begin{equation*}
\log \frac{\log m}{\mu_\epsilon \mutilde_\epsilon^{t_m-1}} = \pp{1-2\frac{\log \mutilde_\epsilon}{\log \mutilde_c}}\log\log m + \log \frac{\mutilde_\epsilon}{\mu_\epsilon} + \theta_m\log\mutilde_\epsilon.
\end{equation*}
Recalling that $\mutilde_\epsilon \to \mutilde_c$ as $\epsilon\to0$ we see that, for sufficiently small $\epsilon$, $\log \mutilde_\epsilon/\log \mutilde_c > 1/2$ and thus $\log m /(\mu_\epsilon \mutilde_\epsilon^{t_m-1}) \to 0$ as $m\to\infty$. It then follows from~\eqref{eq:YlogUpper} that, for such $\epsilon$,
\begin{equation}
\liminf_{m\to\infty} \bbPDo(Y^{(m)}_{t_m} > \log m) \geq \bbP(\calW^\epsilon > 0) = \bbP(\Yhat^\epsilon=\infty).
\label{eq:Yepsliminf}
\end{equation}
Now, $\bc^\epsilon \to \bc$ and $\bctilde^\epsilon \to \bctilde$ as $\epsilon\downarrow0$, so letting $\epsilon\downarrow0$ in~\eqref{eq:Yepsliminf} and using Lemma~\ref{lem:BPTPcgce}(ii) yields
\begin{equation}
\liminf_{m\to\infty} \bbPDo(Y^{(m)}_{t_m} > \log m) \geq \bbP(\Yhat=\infty).
\label{eq:Yliminf}
\end{equation}
Now, for $k=1,2,\ldots$,
\begin{eqnarray*}
\limsup_{m\to\infty} \bbPDo(Y^{(m)}_{t_m} >\log m) & \leq & \limsup_{m\to\infty} \bbPDo(\Yhat^{(m)} >\log m) \\
 & \leq & \limsup_{m\to\infty} \bbPDo(\Yhat^{(m)} >k) \\
 & = & \bbP(\Yhat>k),
\end{eqnarray*}
using Lemmas~\ref{lem:SLLN}(iii) and~\ref{lem:BPTPcgce}(i). Letting $k\to\infty$ then yields
\begin{equation*}
\limsup_{m\to\infty} \bbPDo(Y^{(m)}_{t_m} >\log m) \leq \bbP(\Yhat=\infty),
\end{equation*}
which, together with~\eqref{eq:Yliminf}, establishes the lemma.

In case (ii), a suitable lower bounding branching process is obtained by setting, for $\epsilon<c_{n\dmaxs} \wedge \ctilde_{n\dmaxs -1}$ (where $a \wedge b = \min(a,b)$), $c^\epsilon_i=c_i+\epsilon/(n\dmax)$ ($i=0,1,\ldots,n\dmax -1$), $c^\epsilon_{n\dmaxs}=c_\dmaxs -\epsilon$, $\ctilde^\epsilon_i=\ctilde_i+\epsilon/(n\dmax -1)$ ($i=0,1,\ldots,n\dmax -2$), $\ctilde^\epsilon_{n\dmaxs -1}=\ctilde_{n\dmaxs -1}-\epsilon$, and~\eqref{eq:Yliminf} follows as above.
\end{proof}

For $m=1,2,\ldots$ and $k=0,1,\ldots$, let $Z^{(m)}_k$ denote the number of infectious households in generation $k$ of $E^{(m)}$.
\begin{lem}
\label{lem:limProbYtildeThat}
Let $\beta$ be as in Lemma~\ref{lem:ThatLogm_beta}. Then, for all $\omega_1\in A_1$,
\begin{equation}
\lim_{m\to\infty} \bbPDo\pp{Z^{(m)}_{t_m} > \log m \comma \That^{(m)}_{t_m+1} < (\log m)^\beta} = \bbP(\Yhat=\infty).
\label{eq:ProbYtildeThat}
\end{equation}
\end{lem}
\begin{proof}
Lemmas~\ref{lem:ThatLogm_beta} and~\ref{lem:limProbYt} show that~\eqref{eq:ProbYtildeThat} holds with $Z^{(m)}_{t_m}$ replaced by $Y^{(m)}_{t_m}$. Application of Lemma~\ref{lem:gh}, with $g(m)=(\log m)^\beta$ and $h(m)=2(\log m)^\beta$ then shows that $\lim_{m\to\infty} \bbPDo(Z^{(m)}_{t_m} = Y^{(m)}_{t_m}) = 1$ and the assertion follows.
\end{proof}

\begin{cor}
\label{cor:fwdTotProgCgce}
\begin{enumerate}
\item[(i)] For all $\omega_1\in A_1$, $\ds \lim_{m\to\infty} \bbPDo(\Zhat^{(m)}>\log m) = \bbP(\Yhat=\infty)$;
\item[(ii)] $\ds \lim_{m\to\infty} \bbP(\Zhat^{(m)}>\log m) = \bbP(\Yhat=\infty)$.
\end{enumerate}
\end{cor}
\begin{proof}
Fix $\omega_1\in A_1$. For $k=1,2,\ldots$,
\begin{eqnarray*}
\limsup_{m\to\infty} \bbPDo(\Zhat^{(m)}>\log m) & \leq & \limsup_{m\to\infty} \bbPDo(\Zhat^{(m)}>k) \\
 & = & \bbP(\Yhat>k) \qquad\mbox{(using Theorem~\ref{thm:threshold}(i)),}
\end{eqnarray*}
and letting $k\to\infty$ yields
\begin{equation*}
\limsup_{m\to\infty} \bbPDo(\Zhat^{(m)}>\log m) \leq \bbP(\Yhat=\infty).
\end{equation*}
Also,
\begin{eqnarray*}
\liminf_{m\to\infty} \bbPDo(\Zhat^{(m)}>\log m) & \geq & \liminf_{m\to\infty} \bbPDo(Z^{(m)}_{t_m} > \log m) \\
 & = & \bbP(\Yhat=\infty) \qquad\mbox{(using Lemmas~\ref{lem:ThatLogm_beta} and~\ref{lem:limProbYtildeThat}),}
\end{eqnarray*}
and assertion (i) follows. Assertion (ii) then follows using the dominated convergence theorem, as in the proof of Theorem~\ref{thm:threshold}(ii).
\end{proof}

Note that Theorem~\ref{thm:threshold} and Corollary~\ref{cor:fwdTotProgCgce} imply that if $(h_m)$ is any sequence of real numbers satisfying $h_m\to\infty$ as $m\to\infty$ and $h_m < \log m$ for all $m$, then $\lim_{m\to\infty} \bbPDo(\Zhat^{(m)}\in [h_m,\log m))=0$ for all $\omega_1\in A_1$ and $\lim_{m\to\infty} \bbP(\Zhat^{(m)}\in [h_m,\log m))=0$. Thus, for $m=1,2,\ldots$, it is natural to define a major outbreak as one which infects at least $\log m$ households, i.e.~as one in which the event $\Gbar^{(m)}=\{\omega\in\Omega \,:\, \Zhat^{(m)}(\omega) >\log m\}$ occurs. Let $G^{(m)}=\{\omega\in\Omega \,:\, Z^{(m)}_{t_m} >\log m \comma \That^{(m)}_{t_m+1} <(\log m)^\beta\}$, where $\beta$ is as in Lemma~\ref{lem:ThatLogm_beta}. Clearly $G^{(m)} \subseteq \Gbar^{(m)}$, and Lemma~\ref{lem:ThatLogm_beta} and Corollary~\ref{cor:fwdTotProgCgce} imply $\lim_{m\to\infty} \bbPDo(\Gbar^{(m)} \setminus G^{(m)}) = 0$ for all $\omega_1 \in A_1$ and $\lim_{m\to\infty} \bbP(\Gbar^{(m)} \setminus G^{(m)}) = 0$, so we can take $G^{(m)}$ as our working definition of a major outbreak.

\subsection{Analysis of backward process}
\label{sec:bwdAnalysis}
\subsubsection{Lower bounding branching processes.}
\label{sec:lowerBoundBP}
We now analyse the `backward' process, which describes the generation-wise growth of the susceptibility set (and its neighbours) of a typical individual that is susceptible at time $t_m$ in the forward process, in order to find the asymptotic probability that such an individual is ultimately infected, given that a major outbreak occurs (i.e.~$Z^{(m)}_{t_m} >\log m$ and $\That^{(m)}_{t_m+1} <(\log m)^\beta$, where $\beta$ is as in Lemma~\ref{lem:ThatLogm_beta}). To this end, it is fruitful to have, for all sufficiently small $\epsilon>0$, a branching process ${}_\epsilon X^{(m)}$ which asymptotically bounds $S^{(m)}$ from below until the susceptibility set covers a proportion $\epsilon$ of the households in the population (cf.~Whittle (1955)). In order to do this, we need an almost sure bound, $\etabar(\epsilon)$, for the proportion of households that are neighbours of the susceptibility set when the size (in terms of households) of the susceptibility set is at most $\epsilon m$, which we now obtain.

Suppose that $D$ has infinite support. Recall the definitions of $p_H(\cdot)$ and $\ptilde_H(\cdot)$ from Section~\ref{sec:notAndLimProc}. Let $k_1=\min\{k\,:\,p_H(k)>0\}$ and $\epsilon_0=1-p_H(k_1)-p_H(k_1+1)$. Then, for $\epsilon \in (0,\epsilon_0)$, let $\kappa(\epsilon) = \max \{k \,:\, \sum_{i=k_1}^k p_H(i) \in(0,1-\epsilon)\}$, $\kappa^*(\epsilon) = \max \{k<\kappa(\epsilon) \,:\, p_H(k)>0\}$ and $\eta(\epsilon) = \sum_{i=\kappa^*(\epsilon)}^\infty \ptilde_H(i)$. (The definition of $\kappa^*(\epsilon)$ requires $\kappa(\epsilon)>k_1$, which in turn requires $\epsilon<\epsilon_0$.) Note that $\eta(\epsilon)\downarrow0$ as $\epsilon\downarrow0$. Let $\etabar(\epsilon)=2n\mu_D\eta(\epsilon)$ and, for $m=1,2,\ldots$, let $H^{(m)}_{(1)},H^{(m)}_{(2)},\ldots,H^{(m)}_{(m)}$ be the order statistics of the household degrees $H_1,H_2,\ldots,H_m$.
\begin{lem}
\label{lem:etas}
For any $\omega_1\in A_1$ and $\epsilon\in(0,\epsilon_0)$,
\begin{equation}
\frac{1}{m\mu_{H}^{(m)}(\omega_1)} \sum_{k=m-[\epsilon m]+1}^m H^{(m)}_{(k)}(\omega_1) \leq \eta(\epsilon)
\label{eq:etaBound}
\end{equation}
and
\begin{equation}
\frac{1}{m}\sum_{k=m-[\epsilon m]+1}^m H^{(m)}_{(k)}(\omega_1) \leq \etabar(\epsilon)
\label{eq:etatildeBound}
\end{equation}
for all sufficiently large $m$.
\end{lem}
\begin{proof}
Fix $\omega_1\in A_1$ and note that, for $k=0,1,\ldots$,
\begin{eqnarray}
\lim_{m\to\infty} \frac{1}{m} \sum_{i=1}^m \ind{H_i(\omega_1)=k}
 & = & \lim_{m\to\infty} \frac{1}{m} \sum_{i=1}^m \sum_{\{\bds \,:\, \abs{\bds}=k\}} \ind{\bDs_i(\omega_1)=\bds} \nonumber \\
 & = & \lim_{m\to\infty} \sum_{\{\bds \,:\, \abs{\bds}=k\}} p_{\bds}^{(m)}(\omega_1) \nonumber \\
 & = & \sum_{\{\bds \,:\, \abs{\bds}=k\}} p_{\bds} \qquad \mbox{(using Lemma~\ref{lem:SLLN}(ii))} \nonumber \\
 & = & p_H(k), \label{eq:pHk}
\end{eqnarray}
whence
\begin{eqnarray*}
\lim_{m\to\infty} \frac{1}{m} \sum_{i=1}^m \ind{H_i(\omega_1) \geq \kappa(\epsilon)+1}
 & = & 1-\lim_{m\to\infty} \frac{1}{m} \sum_{i=1}^m \sum_{j=0}^{\kappa(\epsilon)} \ind{H_i(\omega_1) =j} \\
 & = & \sum_{j=\kappa(\epsilon)+1}^\infty p_H(j).
\end{eqnarray*}
Thus, since $\sum_{j=\kappa(\epsilon)+1}^\infty p_H(j) >\epsilon$ (by the definition of $\kappa(\epsilon)$), we have, for all sufficiently large $m$, say $m\geq N_0(\epsilon,\omega_1)$, that $m^{-1}\sum_{i=1}^m \ind{H_i(\omega_1) \geq \kappa(\epsilon)+1}>\epsilon$, whence $H_{(m-[\epsilon m]+1)}^{(m)}(\omega_1) > \kappa(\epsilon)$. Hence, for $m\geq N_0(\epsilon,\omega_1)$,
\begin{eqnarray*}
\frac{1}{m\mu_{H}^{(m)}(\omega_1)} \sum_{k=m-[\epsilon m]+1}^m H^{(m)}_{(k)}(\omega_1)
 & \leq & \frac{1}{m\mu_{H}^{(m)}(\omega_1)} \sum_{i=1}^m \sum_{k=\kappa(\epsilon)+1}^\infty k \ind{H_i(\omega_1)=k} \\
 & = & \frac{m\mu_{H}^{(m)}(\omega_1) - \sum_{i=1}^m \sum_{k=0}^{\kappa(\epsilon)} k \ind{H_i(\omega_1)=k}}{m\mu_{H}^{(m)}(\omega_1)} \\
 & = & 1- \frac{1}{\mu_{H}^{(m)}(\omega_1)} \sum_{k=0}^{\kappa(\epsilon)} \frac{k}{m} \sum_{i=1}^m \ind{H_i(\omega_1)=k} \\
 & \to & 1- \sum_{k=1}^{\kappa(\epsilon)} \ptilde_H(k) = \sum_{k=\kappa(\epsilon)+1}^\infty \ptilde_H(k)
\end{eqnarray*}
as $m\to\infty$, using~\eqref{eq:pHk} and Lemma~\ref{lem:SLLN}(i). Assertion~\eqref{eq:etaBound} follows upon recalling the definition of $\kappa^*(\epsilon)$. The second assertion~\eqref{eq:etatildeBound} follows from the first assertion after applying Lemma~\ref{lem:SLLN}(i) and recalling the definition of $\etabar(\epsilon)$.
\end{proof}

\begin{remsnn}
\begin{enumerate}
\item If $\dmax<\infty$ (i.e.~$D$ has finite support) then it is readily seen that Lemma~\ref{lem:etas} holds with $\eta(\epsilon)=2\dmax\epsilon/\mu_D$ and $\etabar(\epsilon)=n\dmax\epsilon$.
\item For $\omega_1\in A_1$, Lemma~\ref{lem:etas} provides, for all sufficiently large $m$, a bound for the number of half-edges that emanate from households in a susceptibility set (and hence also for the number of households neighbouring a susceptibility set), if the susceptibility set contains no more than $\epsilon m$ households. The number of such half-edges, $H^{(m)}(\epsilon)$ say, is given by the sum of the degrees of the households in the susceptibility set, which is bounded by the sum of the degrees of the $[\epsilon m]$ households of highest degree. Thus, by~\eqref{eq:etatildeBound}, $H^{(m)}(\epsilon) \leq m\etabar(\epsilon)$ for all sufficiently large $m$.
\end{enumerate}
\end{remsnn}

Recall from Section~\ref{sec:bwdConstruction} that the coupling of the susceptibility set process $S^{(m)}$ and its approximating branching process $X^{(m)}$ breaks down when a half-edge is sampled that emanates from an appropriate `bad' set of households. This can happen in two fundamentally different ways. First, a half-edge through which we try to extend the susceptibility set may be paired up with another half-edge through which we want to extend the susceptibility set in the same generation. Note that in this case, neither of the two half-edges concerned actually extends the susceptibility set. Second, the half-edge may be paired with a bad half-edge which is not one through which we wish to extend the susceptibility set in the current generation, in which case the susceptibility set may still be extended, though the offspring distribution is different to that in the branching process. We treat these two cases sequentially.

For $m=1,2,\ldots$ and $k=0,1,\ldots$, let $\Xhat^{(m)}_k=\sum_{i=0}^k X^{(m)}_i$ be the total number of individuals that have lived in the approximating branching process $X^{(m)}$ by time $k$ and let $\What^{(m)}_k=\sum_{i=0}^k S^{(m)}_i$ be the total number of households in the susceptibility set process $S^{(m)}$ up to and including generation $k$. Further, let $\Xhat^{(m)}=\sum_{i=0}^\infty X^{(m)}_i$ and $\What^{(m)}=\sum_{i=0}^\infty S^{(m)}_i$. Suppose that $\omega_1\in A_1$. Then, for all sufficiently large $m$, while $\What^{(m)}_k \leq \epsilon m$, the probability that a half-edge is paired with another half-edge through which we want to extend the susceptibility set in the same generation is no more than $\eta(\epsilon)$. For such $m$, suppose that at some generation $k$ there are $X^{(m)}_{k-1}=i$ `live' half-edges through which we attempt to extend the susceptibility set. Denote by $Y_L$ the number of these half-edges that do not pair up with another of these $i$ live half-edges and let $\Ycheck_L \sim \mbox{Bin}(i,1-\sqrt{\eta(\epsilon)})$. We now show that $Y_L \geqst \Ycheck_L$.

First, define another random variable $\Yhat_L$ as follows. Take a live half-edge, then with probability $\eta(\epsilon)$ pair it up with another live half-edge, otherwise it `survives' to be connected with a non-live half-edge. Repeat this process until all live half-edges have been either paired up or designated to survive. Note that if there is a single live half-edge left at the end of this procedure, it must survive. Let $\Yhat_L$ be the number of surviving half-edges under this regime. Since the proportion of half-edges that are actually live is less than $\eta(\epsilon)$, $Y_L \geqst \Yhat_L$. We now show that $\Yhat_L \geqst \Ycheck_L$ by describing these two random variables as the number of renewals of a discrete-time renewal process by time $i$ and showing that the corresponding lifetime distributions, $\That$ and $\Tcheck$ say, satisfy $\That \leqst \Tcheck$. This we achieve by taking a lifetime in the renewal process as being the number of half-edges examined to find a surviving half-edge.
It is immediate that $\bbP(\Tcheck=k)=(1-\eta(\epsilon)^{\frac12}) \eta(\epsilon)^{\frac{k-1}{2}}$, $k=1,2,\ldots$. Now, since pairing one live half-edge with another obviously uses up two half-edges, $\That$ cannot take even values and $\bbP(\That=2k+1)=(1-\eta(\epsilon)) \eta(\epsilon)^k$, $k=0,1,\ldots$. Elementary calculation shows that $\bbP(\That\geq k)\leq \bbP(\Tcheck\geq k)$, $k=1,2,\ldots$, so $\Tcheck \geqst \That$, whence $\Yhat_L \geqst \Ycheck_L$.

The above argument shows that, in a given generation, the number of half-edges that survive to be paired with non-live half-edges is stochastically larger than if they survive independently with probability $1-\sqrt{\eta(\epsilon)}$. Now consider a live half-edge that survives this first stage and thus is paired with a half-edge chosen uniformly at random from all the non-live half-edges. The probability that it avoids being paired with a half-edge from the bad set is therefore larger than if it were paired with a half-edge chosen uniformly at random from all of the half-edges. Recall that $\ptilde^{(m)}_{\bds}$ is the probability that a half-edge chosen at random from all $m\mu_{H}^{(m)}$ half-edges in the population emanates from a household of type $\bd$. Further, for $\omega_1\in A_1$ and $m$ sufficiently large, conditional on choosing a household of type $\bd$, if $\What^{(m)}_k \leq \epsilon m$ and $Y^{(m)}_{t_m}<(\log m)^\beta$ then the probability of choosing a bad household is bounded above by
\begin{equation*}
\gammatilde_{\bds}^{(m)}(\epsilon) = \frac{(\log m)^\beta + \epsilon m + \etabar(\epsilon)m}{m\ptilde^{(m)}_{\bds}} \wedge 1.
\end{equation*}
This bound is obtained by noting that, under the stated conditions, there are fewer than $(\log m)^\beta$ bad households from the forward process, fewer than $\epsilon m$ households in the susceptibility set and fewer than $\etabar(\epsilon)m$ households that are neighbours of the susceptibility set; and then assuming that all of these bad households are of type $\bd$.

It follows from this discussion that, for $m$ sufficiently large and if $Y^{(m)}_{t_m}<(\log m)^\beta$, then while $\What^{(m)}_k \leq \epsilon m$ the susceptibility set process $S^{(m)}$ is stochastically larger than a branching process, ${}_\epsilon X^{(m)}$ say, in which each potential birth (live half-edge) is aborted independently with probability $\sqrt{\eta(\epsilon)}$ and the potential offspring (live half-edges) of an unaborted birth are obtained by first sampling $\bd$ according to $\ptilde^{(m)}_{\bds}$, then with probability $\gammatilde_{\bds}^{(m)}(\epsilon)$ this unaborted birth is aborted at this stage and otherwise its potential offspring is distributed according to the random variable $\Psitilde_{\bds}$ defined at the end of the paragraph following~\eqref{eq:PhitildedDef}.

The number, $\Xbar^{(m)}_1$ say, of potential births that emanate from the initial individual in the susceptibility set may be found as follows. First a household is chosen uniformly at random from the households not infected by time $t_m$ in the forward process. Suppose that this household is of type $\bd$. Then, if this household is not a neighbour of a household in the forward process, $\Xbar^{(m)}_1$ is distributed according to the random variable $\Psi_\bds$, also defined in the paragraph immediately following~\eqref{eq:PhitildedDef}. If the sampled household is a neighbour of a household in the forward process then $\Xbar^{(m)}_1$ has a different distribution. Suppose that $\That^{(m)}_{t_m+1} < (\log m)^\beta$. Then the number of households that are neighbours of the forward process is less than $2(\log m)^\beta$ and it follows that $\Xbar^{(m)}_1$ is stochastically larger than a random variable, $\Xbbar^{(m)}_1$ say, obtained by first sampling $\bd$ according to $p^{(m)}_{\bds}$ and then setting $\Xbbar^{(m)}_1=0$ with probability $\gamma_{\bds}^{(m)} = \frac{2(\log m)^\beta}{mp^{(m)}_{\bdss}} \wedge 1$, otherwise $\Xbbar^{(m)}_1$ is distributed according to $\Psi_{\bds}$.

Assume that there is a single ancestor in the branching process ${}_\epsilon X^{(m)}$, which has a number of potential offspring distributed as $\Xbbar^{(m)}_1$. We now have a complete description of how ${}_\epsilon X^{(m)}$ evolves. Let ${}_\epsilon \Xhat^{(m)}$ and ${}_\epsilon \What^{(m)}$ be, respectively, the total number of potential and unaborted births in ${}_\epsilon X^{(m)}$. Recall the event $G^{(m)}$ defined at the end of Section~\ref{sec:earlyMajOB}, giving our working definition of a major outbreak. The above arguments show that
\begin{eqnarray}
\bbPDo(\What^{(m)}>[\epsilon m] \cond G^{(m)}) & \geq & \bbPDo({}_\epsilon \What^{(m)} \geq [\epsilon m]) \nonumber \\
 & \geq & \bbPDo({}_\epsilon \What^{(m)} =\infty) \nonumber \\
 & = & \bbPDo({}_\epsilon \Xhat^{(m)} =\infty). \label{eq:W}
\end{eqnarray}

For the branching process ${}_\epsilon X^{(m)}$, let ${}_\epsilon \bb^{(m)}=({}_\epsilon b^{(m)}_0,{}_\epsilon b^{(m)}_1,\ldots)$ denote the distribution of the number of potential offspring of the initial individual and let ${}_\epsilon \bbtilde^{(m)}=({}_\epsilon \btilde^{(m)}_0,{}_\epsilon \btilde^{(m)}_1,\ldots)$ denote the distribution of the number of potential offspring of a typical potential birth. Then
\begin{eqnarray*}
{}_\epsilon b^{(m)}_0 & = & \sum_{\bds\in\bbZ_+^n} p_{\bds}^{(m)} \pp{\gamma_{\bds}^{(m)} + (1-\gamma_{\bds}^{(m)})\bbP(\Psi_{\bds}=0)}, \\
{}_\epsilon b^{(m)}_k & = & \sum_{\bds\in\bbZ_+^n} p_{\bds}^{(m)} (1-\gamma_{\bds}^{(m)})\bbP(\Psi_{\bds}=k) \qquad (k=1,2,\ldots), \\
{}_\epsilon \btilde^{(m)}_0 & = & \sqrt{\eta(\epsilon)} + (1-\sqrt{\eta(\epsilon)}) \sum_{\bds\in\bbZ_+^n} \ptilde_{\bds}^{(m)} \pp{\gammatilde_{\bds}^{(m)}(\epsilon) + (1-\gammatilde_{\bds}^{(m)}(\epsilon))\bbP(\Psitilde_{\bds}=0)}
\end{eqnarray*}
and
\begin{equation*}
{}_\epsilon \btilde^{(m)}_k = (1-\sqrt{\eta(\epsilon)}) \sum_{\bds\in\bbZ_+^n} \ptilde_{\bds}^{(m)} (1-\gammatilde_{\bds}^{(m)}(\epsilon))\bbP(\Psitilde_{\bds}=k) \qquad (k=1,2,\ldots).
\end{equation*}
Note that ${}_\epsilon \bb^{(m)}$ does not depend on $\epsilon$, however it is distinct from $\bb^{(m)}$ and we retain the notation ${}_\epsilon \bb^{(m)}$ to indicate that it is associated with the branching process ${}_\epsilon X^{(m)}$.

The following lemma is useful for determining the limits of the distributions ${}_\epsilon \bb^{(m)}$ and ${}_\epsilon \bbtilde^{(m)}$ as $m\to\infty$. Its proof is standard and is hence omitted.

\begin{lem}
\label{lem:sc}
Suppose that, for all $\bd\in\bbZ_+^n$ and $m=1,2,\ldots$, the real numbers
\begin{enumerate}
\item[(i)] $p_{\bds}^{(m)}$ and $p_{\bds}$ are non-negative and satisfy $p_{\bds}^{(m)} \to p_{\bds}$ as $m\to\infty$ and $\sum_{\bds\in\bbZ_+^n} p_{\bds}^{(m)} = \sum_{\bds\in\bbZ_+^n} p_{\bds} =1$;
\item[(ii)] $\alpha^{(m)}_{\bds}$ and $\alpha_{\bds}$ belong to $[0,1]$ and satisfy $\alpha^{(m)}_{\bds} \to \alpha_{\bds}$ as $m\to\infty$;
\item[(iii)] $c_\bds$ belong to $[0,1]$.
\end{enumerate}
Then, as $m\to\infty$,
\begin{equation*}
\sum_{\bds\in\bbZ_+^n} p_{\bds}^{(m)} \alpha^{(m)}_{\bds} c_\bds \to \sum_{\bds\in\bbZ_+^n} p_{\bds} \alpha_{\bds} c_\bds.
\end{equation*}
\end{lem}

For $\bd\in\bbZ_+^n$ and $\epsilon\in(0,\epsilon_0)$, let
\begin{equation*}
\gammatilde_\bds(\epsilon) = \begin{cases}
\pp{\frac{\epsilon+\etabar(\epsilon)}{\ptilde_\bdss}} \wedge 1 & \mbox{if $\ptilde_{\bds}>0$, }\\
0 & \mbox{if $\ptilde_\bds=0$.}\end{cases}
\end{equation*}
\begin{lem}
\label{lem:epstoZeroOffspProbs}
For all $\omega_1\in A_1$, $\lim_{m\to\infty} {}_\epsilon \bb^{(m)} = \bb$ and $\lim_{m\to\infty} {}_\epsilon \bbtilde^{(m)} = {}_\epsilon \bbtilde$, where $\bb=(b_0,b_1,\ldots)$ is as in Section~\ref{sec:notAndLimProc}, and ${}_\epsilon \bbtilde=({}_\epsilon \btilde_0,{}_\epsilon \btilde_1,\ldots)$ is given by
\begin{equation*}
{}_\epsilon \btilde_0 = \sqrt{\eta(\epsilon)} + (1-\sqrt{\eta(\epsilon)}) \sum_{\bds\in\bbZ_+^n} \ptilde_{\bds} \pp{\gammatilde_{\bds}(\epsilon) + (1-\gammatilde_{\bds}(\epsilon))\bbP(\Psitilde_{\bds}=0)}
\end{equation*}
and
\begin{equation*}
{}_\epsilon \btilde_k = (1-\sqrt{\eta(\epsilon)}) \sum_{\bds\in\bbZ_+^n} \ptilde_{\bds} (1-\gammatilde_{\bds}(\epsilon))\bbP(\Psitilde_{\bds}=k) \qquad (k=1,2,\ldots).
\end{equation*}
\end{lem}
\begin{proof}
Note that, for $\omega_1\in A_1$, $\gamma^{(m)}_\bds(\omega_1) \to0$ and $\gammatilde^{(m)}_\bds(\epsilon,\omega_1) \to \gamma_\bds(\epsilon)$ as $m\to\infty$ (for all $\bd$ with $p_\bds>0$). The required assertions then follow using Lemma~\ref{lem:sc}.
\end{proof}
\begin{remnn}
It is easily verified that $\sum_{k=0}^\infty {}_\epsilon \btilde_k = 1$, i.e.~that ${}_\epsilon \bbtilde$ is a proper probability distribution.
\end{remnn}

Recall the definition of $\epsilon_0$ in the paragraph preceding Lemma~\ref{lem:etas} and, for $\epsilon\in(0,\epsilon_0)$, let ${}_\epsilon X=({}_\epsilon X_k,\, k=0,1,\ldots)\sim {\rm BP}(1,\bb,{}_\epsilon \bbtilde)$.  Let ${}_\epsilon \Xhat$ denote the total progeny of ${}_\epsilon X$, excluding the ancestor. Let $(\Xhat,\Xhat_A)$ denote the total progeny of the branching process $(X,X_A)$ (defined at the end of Section~\ref{sec:notAndLimProc}), including the ancestor. Also let $\Xhat_A^{(m)} = \sum_{i=0}^\infty X_{Ai}^{(m)}$, so $(\Xhat^{(m)},\Xhat_A^{(m)})$ is the total progeny of $(X^{(m)},X_A^{(m)})$.

\begin{lem}
\label{lem:bwdTotProgCgce}
\begin{enumerate}\item[(i)] For all $\omega_1\in A_1$,
\begin{enumerate}
\item[(a)] $\ds \lim_{m\to\infty} \bbPDo(\Xhat^{(m)} + \Xhat_A^{(m)} =k) = \bbP(\Xhat + \Xhat_A=k) \qquad (k=1,2,\ldots)$;
\item[(b)] $\ds \lim_{m\to\infty} \bbPDo(\Xhat^{(m)} + \Xhat_A^{(m)} =\infty) = \bbP(\Xhat=\infty)$.
\end{enumerate}
\item[(ii)] For all $\omega_1\in A_1$ and $\epsilon\in(0,\epsilon_0)$,
\begin{enumerate}
\item[(a)] $\ds \lim_{m\to\infty} \bbPDo({}_\epsilon \Xhat^{(m)}=k) = \bbP({}_\epsilon \Xhat=k) \qquad (k=1,2,\ldots)$;
\item[(b)] $\ds \lim_{m\to\infty} \bbPDo({}_\epsilon \Xhat^{(m)}=\infty) = \bbP({}_\epsilon \Xhat=\infty)$.
\end{enumerate}
\end{enumerate}
\end{lem}
\begin{proof}
For all $\omega_1\in A_1$ and $\bd\in\bbZ_+^n$, $p^{(m)}_\bds(\omega_1) \to p_\bds$ and $\ptilde^{(m)}_\bds(\omega_1) \to \ptilde_\bds$ as $m\to\infty$, so, using \Scheffe's theorem, $\bb^{(m)}(\omega_1) \to \bb$ and $\bbtilde^{(m)}(\omega_1) \to \bbtilde$ as $m\to\infty$. Part~(ii)(b) then follows using Lemma~\ref{lem:BPTPcgce}(ii) and noting that, almost surely, $\Xhat^{(m)} + \Xhat_A^{(m)} =\infty$ if and only if $\Xhat^{(m)} =\infty$. A similar argument shows that, for all $\omega_1\in A_1$, the offspring laws of $(X^{(m)},X_A^{(m)})$ converge to those of $(X,X_A)$ as $m\to\infty$. Part~(i)(a) then follows from the extension of Lemma~\ref{lem:BPTPcgce}(i) to two-type branching processes. Part~(ii) of the lemma follows immediately from Lemmas~\ref{lem:epstoZeroOffspProbs} and~\ref{lem:BPTPcgce}.
\end{proof}

\subsubsection{Relative final size of a major outbreak.}
\label{sec:majOBInfProp}
For $m=1,2,\ldots$, let $B^{(m)}$ be the event that an individual chosen uniformly at random from all individuals that are susceptible at time $t_m$ in the forward process is ultimately infected by the epidemic $E^{(m)}$. Thus, if $\calA^{(m)}$ denotes the set of global neighbours of $\calS^{(m)}$ then $B^{(m)}$ occurs if and only if one of the $Z^{(m)}_{t_m}$ `live' half-edges from the forward process is paired in the construction of $\calS^{(m)} \cup \calA^{(m)}$. Recall the working definition of a major outbreak, viz. $G^{(m)} = \{ Z^{(m)}_{t_m}>\log m \comma \That^{(m)}_{t_m+1}<(\log m)^\beta \}$, where $\beta$ is as in Lemma~\ref{lem:ThatLogm_beta}.
\begin{thm}
\label{thm:majOBInfProb}
For all $\omega_1\in A_1$,
\begin{equation*}
\lim_{m\to\infty} \bbPDo(B^{(m)} \cond G^{(m)}) = \bbP(\Xhat=\infty).
\end{equation*}
\end{thm}
\begin{proof}
For $m=1,2,\ldots$, let $T_P^{(m)}$ be the number of half-edge pairings made in the construction of $\calS^{(m)} \cup \calA^{(m)}$ until one of the $Z^{(m)}_{t_m}$ live half-edges from the forward process is chosen. In determining $T_P^{(m)}$ it is assumed that, if necessary, the pairings continue after $\calS^{(m)} \cup \calA^{(m)}$ goes extinct and that $T_P^{(m)}$ includes the pairing when the first live half-edge is chosen.

Fix $\omega_1\in A_1$. First we obtain an upper bound for $\bbPDo(B^{(m)} \cond G^{(m)})$. For all fixed $k\in\bbN$,
\begin{equation}
1-\bbPDo(B^{(m)} \cond G^{(m)}) \geq \bbPDo(T_P^{(m)} >k \comma \Xhat^{(m)}+\Xhat_A^{(m)}\leq k \comma \taubar^{(m)}>k \cond G^{(m)}),
\label{eq:1-PBcondGLower}
\end{equation}
where $\taubar^{(m)}$ is the number of households in the construction of $\calS^{(m)} \cup \calA^{(m)}$ when the first bad half-edge is chosen. Note that $\taubar^{(m)}=1$ if the initial individual in $(X^{(m)},X_A^{(m)})$ belongs to the set of bad households at time $t_m$ in the forward process. Given $G^{(m)}$, the number of such bad households is less than $(\log m)^\beta$, so $\bbPDo(\taubar^{(m)}=1 \cond G^{(m)})\to0$ as $m\to\infty$. Arguing as in the proof of Theorem~\ref{thm:threshold} then shows that, for all $k\in\bbN$,
\begin{equation}
\lim_{m\to\infty} \bbPDo(\taubar^{(m)} >k \cond G^{(m)}) = 1.
\label{eq:taubarToInfty}
\end{equation}

Let $Q^{(m)}$ denote the number of half-edges used up to time $t_m$ in the forward process. Now, for all $k\in\bbN$,
\begin{equation}
\bbPDo(T_P^{(m)} >k \cond G^{(m)} \comma Q^{(m)} \comma Z^{(m)}_{t_m}) = \prod_{i=1}^k \pp{ \frac{m\mu_h^{(m)}(\omega_1) - Q^{(m)} -2(i-1) - Z^{(m)}_{t_m}}{m\mu_h^{(m)}(\omega_1) - Q^{(m)} -2(i-1)} }
\label{eq:TPcalcn}
\end{equation}
and, since we have conditioned on $G^{(m)}$, $Q^{(m)}<2(\log m)^\beta$ and $Z^{(m)}_{t_m}<2(\log m)^\beta$. It then follows from~\eqref{eq:TPcalcn} that
\begin{equation}
\lim_{m\to\infty} \bbPDo(T_P^{(m)} >k \cond G^{(m)} ) = 1
\label{eq:TPtoInfty}
\end{equation}
for all $k\in\bbN$. Letting $m\to\infty$ in~\eqref{eq:1-PBcondGLower}, using~\eqref{eq:taubarToInfty} and~\eqref{eq:TPtoInfty}, and noting that $\Xhat^{(m)}+\Xhat_A^{(m)}$ and $G^{(m)}$ are conditionally independent given $\bD(\omega_1)$, yields, for all $k\in\bbN$,
\begin{eqnarray*}
\limsup_{m\to\infty} \bbPDo(B^{(m)} \cond G^{(m)}) & \leq & \limsup_{m\to\infty} \bbPDo(\Xhat^{(m)}+\Xhat_A^{(m)}>k) \\
 & = & \bbP(\Xhat+\Xhat_A>k),
\end{eqnarray*}
using Lemma~\ref{lem:bwdTotProgCgce}(i)(a). Letting $k\to\infty$ then yields
\begin{equation}
\limsup_{m\to\infty} \bbPDo(B^{(m)} \cond G^{(m)}) \leq \bbP(\Xhat=\infty).
\label{eq:limsupBound}
\end{equation}

Now we obtain a lower bound for $\bbPDo(B^{(m)} \cond G^{(m)})$. First note, using~\eqref{eq:TPcalcn}, that for any $\epsilon\in(0,1)$, we have
\begin{equation*}
\bbPDo(T_P^{(m)}>[\epsilon m] \cond G^{(m)}) \leq \pp{ 1 - \frac{\log m}{m\mu_H^{(m)}(\omega_1)} }^{[\epsilon m]} \leq \exp \pp{ \frac{-[\epsilon m]\log m}{m\mu_H^{(m)}(\omega_1)} }.
\end{equation*}
Now, $\mu_H^{(m)}(\omega_1) \to n\mu_D$ as $m\to\infty$ (since $\omega_1\in A_1$), so $[\epsilon m]\log m/m\mu_H^{(m)}(\omega_1) \to \infty$ as $m\to\infty$, whence
\begin{equation}
\lim_{m\to\infty} \bbPDo(T_P^{(m)} \leq [\epsilon m] \cond G^{(m)}) = 1.
\label{eq:TPleqEpsm}
\end{equation}
Also note that, since $\calS^{(m)}$ is obviously contained in $\calS^{(m)} \cup \calA^{(m)}$,
\begin{equation}
\bbPDo(B^{(m)} \cond G^{(m)}) \geq \bbPDo(T_P^{(m)} \leq [\epsilon m] \comma \What^{(m)}>[\epsilon m] \cond G^{(m)}),
\label{eq:PBcondGLower}
\end{equation}
for any $\epsilon \in (0,1)$. Thus, using~\eqref{eq:PBcondGLower} and~\eqref{eq:TPleqEpsm}, then~\eqref{eq:W} and Lemma~\ref{lem:bwdTotProgCgce}(ii)(b), for any $\epsilon\in(0,\epsilon_0)$,
\begin{eqnarray}
\liminf_{m\to\infty} \bbPDo(B^{(m)} \cond G^{(m)}) & \geq & \liminf_{m\to\infty} \bbPDo(\What^{(m)}>[\epsilon m] \cond G^{(m)}) \nonumber \\
 & \geq & \liminf_{m\to\infty} \bbPDo({}_\epsilon\Xhat^{(m)} =\infty) \nonumber \\
 & = & \bbP({}_\epsilon\Xhat =\infty). \label{eq:epsLiminfBound}
\end{eqnarray}
It is easily verified, using the dominated convergence theorem, that $({}_\epsilon \bb,{}_\epsilon \bbtilde) \to (\bb,\bbtilde)$ as $\epsilon\downarrow0$, so letting $\epsilon\downarrow0$ in~\eqref{eq:epsLiminfBound} and using Lemma~\ref{lem:BPTPcgce}(ii) yields
\begin{equation*}
\liminf_{m\to\infty} \bbPDo(B^{(m)} \cond G^{(m)}) \geq \bbP(\Xhat =\infty),
\end{equation*}
which together with~\eqref{eq:limsupBound} establishes the assertion of the theorem.
\end{proof}

For $m=1,2,\ldots$, let $\Zbar^{(m)}_k$ be the total number of individuals infected by time $k$ in the forward epidemic process $E^{(m)}$ ($k=0,1,\ldots$) and let $\Zbar^{(m)}$ denote the total number of individuals who are ultimately infected in $E^{(m)}$.
\begin{cor}
\label{cor:UncondMajOBInfProb}
\begin{enumerate}
\item[(i)] For all $\omega_1\in A_1$, $\ds \lim_{m\to\infty} \frac{1}{mn} \bbEDo[\Zbar^{(m)} \cond G^{(m)}] = \bbP(\Xhat=\infty)$;
\item[(ii)] $\ds \lim_{m\to\infty} \frac{1}{mn} \bbE[\Zbar^{(m)} \cond G^{(m)}] = \bbP(\Xhat=\infty)$.
\end{enumerate}
\end{cor}
\begin{proof}
Fix $\omega_1\in A_1$. For $m=1,2,\ldots$, let $\Xbar_{t_m}$ denote the number of susceptible individuals at time $t_m$ in the forward process, and label these individuals $1,2,\ldots,\Xbar_{t_m}$. Then
\begin{equation*}
\Zbar^{(m)} = \Zbar^{(m)}_{t_m} + \sum_{i=1}^{\Xbar_{t_m}} \ind{\mbox{\scriptsize $i$ ultimately infected}}.
\end{equation*}
Given the occurrence of $G^{(m)}$, $\Zbar^{(m)}_{t_m} < 2n(\log m)^\beta$ and $\Xbar_{t_m} > nm - 2n(\log m)^\beta$. Thus
\begin{equation*}
\lim_{m\to\infty} \frac{1}{mn} \bbEDo[\Zbar^{(m)} \cond G^{(m)}] = \lim_{m\to\infty} \bbPDo(B^{(m)} \cond G^{(m)})
\end{equation*}
and assertion (i) follows using Theorem~\ref{thm:majOBInfProb}. Assertion (ii) then follows by the dominated convergence theorem.
\end{proof}

Finally, note from the discussion at the end of Section~\ref{sec:earlyMajOB} that Corollary~\ref{cor:UncondMajOBInfProb} holds with $G^{(m)}$ replaced by $\Gbar^{(m)}$, where $\Gbar^{(m)}$ is the event that the epidemic $E^{(m)}$ infects at least $\log m$ households.

\section{Concluding comments}
\label{sec:Discussion}
We have analysed the spread of an SIR epidemic within a population structure that features some significant departures from traditional homogeneous mixing; specifying both a local household structure and using random networks with an arbitrary degree distribution (with finite variance) to model potential `global' contacts. Rigorous limit theorems were obtained, valid as the number of households $m\to\infty$, from which one can determine the probability of a major outbreak and the expected relative final size of such an outbreak. The potential usefulness of these results was verified by showing, numerically, that these asymptotic results provide good approximations for the behaviour of moderately sized finite populations.

As stated in Section~\ref{sec:Model}, our results easily generalise to allow for unequal household sizes. For example, we can decompose $R_*$ in a variable household size  framework as $R_* = \sum_{n=1}^\infty \rhotilde_n R_*^{(n)}$, where $\rhotilde_n$ is the size-biased proportion of households of size $n$ and $R_*^{(n)}$ is the threshold parameter $R_*$ in the case of a fixed household size $n$. (The size-bias of $\rhotilde_n$ arises because if a proportion $\rho_n$ of households are of size $n$ then an individual chosen uniformly at random is in a household of size $n$ with probability proportional to $n\rho_n$; thus we require $\sum_{n=1}^\infty n\rho_n <\infty$.) Full details of this generalisation will appear in a forthcoming paper, which will discuss our model from a more applied viewpoint.

Another condition that we have required is that the variance, $\sigma^2_D$, of the degree distribution is finite. Whilst this is necessary for all of our proofs, the PGFs of $C$, $\Ctilde$, $B$ and $\Btilde$ are all well-defined so long as $\mu_D<\infty$ and numerical studies (along the lines of those encompassed by Figure~\ref{fig:cgce}) indicate that our methods at least give good approximations when $\sigma^2_D=\infty$. This is particularly relevant in light of several of the studies cited by Newman (2003, Section~III.C), which suggest that degree distributions which asymptotically follow some power law are appropriate models in some real-world situations. We note, however, that when $\sigma^2_D=\infty$ it is not known (to our knowledge) whether self-loops and parallel edges remain sufficiently sparse in the network, so the argument that our results continue to hold if we condition on there being no such imperfections (second paragraph of Section~\ref{sec:Model}) may not be valid.

Of course there are other features of our model that in many circumstances will be unrealistic. In particular, the method of construction of the random graph---pairing the half-edges uniformly at random---ensures not only that there are (asymptotically) very few 1-cycles (self-loops) and 2-cycles (parallel edges) in the resulting multigraph, but also that there are very few 3-cycles (triangles). Thus, in the asymptotic model that we analyse, individuals have no mutual acquaintances outside their household, which is unrealistic. Similarly the random graph model has very few edges which join individuals in the same pair of households, i.e.\ the acquaintances of two individuals are, with probability close to 1, all in distinct households. That this is the case stems from the construction of the random graph: although there is heterogeneity amongst the individuals (through differing degrees), the uniformly at random pairing of half-edges means that the mixing is still homogeneous---this being critical for the branching process approximations. In this sense it seems fair to say that our model incorporates some heterogeneity of both the individuals in the population (via the differing degrees of individuals and varying household sizes) and their mixing (having both local and global infection).

Nevertheless, our model does capture some important heterogeneities which are present in real populations and which doubtless have a significant effect on the spread of disease through these populations. Some additional features, such as having the degree distribution $D$ or the infection rates $\lambda_L$ and $\lambda_G$ depend on household size or incorporating correlation between the degrees of individuals within the same household can in principle be included in our model relatively simply, though the calculations quickly become very cumbersome.

The usual approach for obtaining fully rigorous results concerning the final size of a major epidemic on a random network is via the existence and uniqueness of a giant component in an associated bond percolation model (see e.g.\ Britton \etal\ (2007) and the discussion in Section 4 of Britton \etal\ (2008)).  This requires that the infectious period is constant (though see Kenah and Robins (2007)) and fully rigorous results concerning the component structure of the percolation model, which may not be easy to prove. We have developed a different approach, which does not require a constant infectious period.  Although not the focus of the paper, it seems plausible that our methods can be used to prove existence and uniqueness of a giant component for our random network (and indeed for other network models) and that they might also be applicable to epidemics on other random graph models, such as the random intersection graph considered by Britton \etal\ (2008).

Further study of this model will include an analysis of the effect of vaccination on epidemic spread (work ongoing) and it seems likely that a central limit theorem for the final size of a major outbreak might be derived using methods similar to those of Ball and Neal (2008).

\paragraph*{Acknowledgements}
This research was supported by the UK Engineering and Physical Sciences Research Council, under research grant number EP/E038670/1 (Frank Ball and David Sirl) and by the Netherlands Organisation for Scientific Research (NWO) through a VICI grant awarded to Ronald Meester (Pieter Trapman).

\section*{References}
\noindent

{\sc Andersson, H.} (1997).
\newblock Epidemics in a population with social structures.
\newblock {\em Math. Biosci.\/} {\bf 140,} 79--84.

{\sc Andersson, H.} (1998).
\newblock Limit theorems for a random graph epidemic model.
\newblock {\em Ann. Appl. Probab.\/} {\bf 8,} 1331--1349.

{\sc Andersson, H.} (1999).
\newblock Epidemic models and social networks.
\newblock {\em Math. Sci.\/} {\bf 24,} 128--147.

{\sc Ball, F.~G.} (1986).
\newblock A unified approach to the distribution of total size and total area
  under the trajectory of infectives in epidemic models.
\newblock {\em Adv. in Appl. Probab.\/} {\bf 18,} 289--310.

{\sc Ball, F.~G. and Lyne, O.~D.} (2001).
\newblock Stochastic multitype {SIR} epidemics among a population partitioned
  into households.
\newblock {\em Adv. in Appl. Probab.\/} {\bf 33,} 99--123.

{\sc Ball, F.~G., Mollison, D. and Scalia-Tomba, G.} (1997).
\newblock Epidemics with two levels of mixing.
\newblock {\em Ann. Appl. Probab.\/} {\bf 7,} 46--89.

{\sc Ball, F.~G. and Neal, P.~J.} (2002).
\newblock A general model for stochastic {SIR} epidemics with two levels of
  mixing.
\newblock {\em Math. Biosci.\/} {\bf 180,} 73--102.

{\sc Ball, F.~G. and Neal, P.~J.} (2003).
\newblock The great circle epidemic model.
\newblock {\em Stochastic Process. Appl.\/} {\bf 107,} 233--268.

{\sc Ball, F.~G. and Neal, P.~J.} (2008).
\newblock Network epidemic models with two levels of mixing.
\newblock {\em Math. Biosci.\/} {\bf 212,} 69--87.

{\sc Ball, F.~G. and O'Neill, P.~D.} (1999).
\newblock The distribution of general final state random variables for
  stochastic epidemic models.
\newblock {\em J. Appl. Probab.\/} {\bf 36,} 473--491.

{\sc Becker, N.~G. and Dietz, K.} (1995).
\newblock The effect of household distribution on transmission and control of
  highly infectious diseases.
\newblock {\em Math. Biosci.\/} {\bf 127,} 207--219.

{\sc Billingsley, P.} (1968).
\newblock {\em Convergence of probability measures}.
\newblock John Wiley \& Sons Inc., New York.

{\sc Britton, T., Deijfen, M., Lager{\aa}s, A.~N. and Lindholm, M.} (2008).
\newblock Epidemics on random graphs with tunable clustering.
\newblock {\em J. Appl. Probab.\/} {\bf 45,} 743--756.

{\sc Britton, T., Janson, S. and Martin-L{\"o}f, A.} (2007).
\newblock Graphs with specified degree distributions, simple epidemics, and
  local vaccination strategies.
\newblock {\em Adv. in Appl. Probab.\/} {\bf 39,} 922--948.

{\sc Durrett, R.} (2006).
\newblock {\em Random graph dynamics}.
\newblock Cambridge Series in Statistical and Probabilistic Mathematics.
  Cambridge University Press, Cambridge.

{\sc Haccou, P., Jagers, P. and Vatutin, V.} (2005).
\newblock {\em Branching processes: Variation, growth, and extinction of
  populations}.
\newblock Cambridge University Press, Cambridge.

{\sc van~der Hofstad, R., Hooghiemstra, G. and Znamenski, D.} (2007).
\newblock Distances in random graphs with finite mean and infinite variance degrees.
\newblock {\em Electron. J. Probab.\/} {\bf 12,} 703--766.

{\sc Janson, S.} (2009).
\newblock The probability that a random multigraph is simple.
\newblock {\em Combin. Probab. Comput.\/} {\bf 18,} 205--225.

{\sc Kenah, E. and Robins, J.~M.} (2007).
\newblock Second look at the spread of epidemics on networks.
\newblock {\em Phys. Rev. E\/} {\bf 76,} 036113.

{\sc Kiss, I.~Z., Green, D.~M. and Kao, R.~R.} (2006).
\newblock The effect of contact heterogeneity and multiple routes of
  transmission on final epidemic size.
\newblock {\em Math. Biosci.\/} {\bf 203,} 124--136.

{\sc Kuulasmaa, K.} (1982).
\newblock The spatial general epidemic and locally dependent random graphs.
\newblock {\em J. Appl. Probab.\/} {\bf 19,} 745--758.

{\sc Newman, M. E.~J.} (2002).
\newblock Spread of epidemic disease on networks.
\newblock {\em Phys. Rev. E\/} {\bf 66,} 016128.

{\sc Newman, M. E.~J.} (2003).
\newblock The structure and function of complex networks.
\newblock {\em SIAM Rev.\/} {\bf 45,} 167--256 (electronic).

{\sc Pellis, L., Ferguson, N.~M. and Fraser, C.} (2008).
\newblock The relationship between real-time and discrete-generation models of epidemic spread.
\newblock {\em Math. Biosci.\/} {\bf 216,} 63--70.

{\sc Trapman, P.} (2007).
\newblock On analytical approaches to epidemics on networks.
\newblock {\em Theor. Pop. Biol.\/} {\bf 71,} 160--173.

{\sc Watts, D.~J. and Strogatz, S.~H.} (1998).
\newblock Collective dynamics of `small-world' networks.
\newblock {\em Nature\/} {\bf 393,} 440--442.

{\sc Whittle, P.} (1955).
\newblock The outcome of a stochastic epidemic---a note on {B}ailey's paper.
\newblock {\em Biometrika\/} {\bf 42,} 116--122.

\end{document}